\documentclass[a4paper,12pt]{article}
\usepackage[latin1]{inputenc}
\usepackage{amssymb}
\usepackage{graphicx}
\usepackage{amsmath}
\usepackage{calc}
\usepackage{eurosym}
\usepackage{anysize}
\usepackage{amsthm}
\usepackage{enumerate}
\usepackage{semtrans}
\usepackage{eufrak}
\usepackage{mathrsfs}
\usepackage{yfonts}
\usepackage{bbold}
\usepackage{color}
\usepackage{pst-node}
\usepackage{pstricks}
\usepackage[nottoc, notlof, notlot]{tocbibind}
\usepackage{wasysym}
\usepackage{MnSymbol}
\marginsize{2,5cm}{2,5cm}{2,2cm}{2,3cm}

\title{Interacting Urn Models}
\date {\today}
\author{Micka{\"{e}}l Launay}

\definecolor{jaune} {cmyk}{0,0.4,1,0}%
\definecolor{bleu} {cmyk}{1,0,0,0}%

\newcommand{\N}{\mathbb{N}}

\newcommand{\R}{\mathbb{R}}
\newcommand{\Z}{\mathbb{Z}}

\newcommand{\RUPa}{{IUM}}
\newcommand{\drawall}{{all the $U$ urns combined}}
\newcommand{\up}{{Urn Path}}
\newcommand{\ups}{{Urn Paths}}

\begin{document}

\newtheorem{lemme}{Lemma}[section]
\newtheorem{theorem}[lemme]{Theorem}
\newtheorem{prop}[lemme]{Proposition}
\newtheorem{coro}[lemme]{Corollaire}
\newtheorem{conj}[lemme]{Conjecture}
\newtheorem{defi}[lemme]{Definition}
\newtheorem{ex}[lemme]{Exemples}

\maketitle

\begin{abstract}
The aim of this paper is to study the asymptotic behavior of strongly reinforced interacting urns 
with partial memory sharing. The reinforcement mechanism considered is as follows: draw at each step and for each urn a white or black ball from either all the urns combined (with probability $p$) or the urn alone (with probability $1-p$) and add a new ball of the same color to this urn. The probability of drawing a ball of a certain color is proportional to $w_k$ where $k$ is the number of balls of this color. The higher the $p$, the more memory is shared between the urns. The main results can be informally stated as follows: in the exponential case $w_k=\rho^k$, if $p\geq 1/2$ then all the urns draw the same color after a finite time, and if $p<1/2$ then some urns fixate on a unique color and others keep drawing both black and white balls.
\end{abstract}

\section{Introduction}

Assume that we are given $U$ urns, where $U \in \N$.
To each urn give a unique label in $\{1,\ldots,U\}$, 
and refer to the urn labeled by $u$ as the ``urn $u$''.
Assume we are given in addition
a {\em reinforcement weight sequence}
$(w_i)_{i\in\N\cup\{0\}}\in{\R_+}^{\N\cup\{0\}}$, and an {\em interaction parameter} $p\in[0,1]$.

Consider the following model:

- At time $n=0$, all the urns are empty;
\begin{center}
	\pspicture(0,-1)(8.5,2.5)
			\psline[border=0pt](0,2)(0,0)
			\psline[border=0pt](0,0)(1,0)
			\psline[border=0pt](1,0)(1,2)
			
			\psline[border=0pt]{-}(1.5,2)(1.5,0)
			\psline[border=0pt]{-}(1.5,0)(2.5,0)
			\psline[border=0pt]{-}(2.5,0)(2.5,2)

			\psline[border=0pt]{-}(3,2)(3,0)
			\psline[border=0pt]{-}(3,0)(4,0)
			\psline[border=0pt]{-}(4,0)(4,2)
			
			\psdots[dotsize=1.5pt](4.5,1)(5,1)(5.5,1)
			
			\psline[border=0pt]{-}(6,2)(6,0)
			\psline[border=0pt]{-}(6,0)(7,0)
			\psline[border=0pt]{-}(7,0)(7,2)
			
			\psline[border=0pt]{-}(7.5,2)(7.5,0)
			\psline[border=0pt]{-}(7.5,0)(8.5,0)
			\psline[border=0pt]{-}(8.5,0)(8.5,2)
			
			\rput(.5,-.5){{$1$}}
			\rput(2,-.5){{$2$}}
			\rput(3.5,-.5){{$3$}}
			\rput(6.5,-.5){{$U-1$}}
			\rput(8,-.5){{$U$}}
						
	\endpspicture

\noindent
{\footnotesize Figure 1:
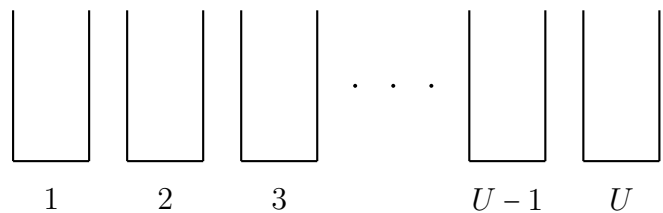
At $n=0$ all the $U$ urns are empty.}
\setcounter{figure}{1}

\vspace{0.3cm}
\noindent

\end{center}

- At time $n$ each of the $U$ urns contains exactly $n$ balls. The balls are either black or white. 
The next figure is an illustration corresponding to time $n=3$:

	\begin{center}
	\pspicture(0,-1)(8.5,2.5)
			\psline[border=0pt](0,2)(0,0)
			\psline[border=0pt](0,0)(1,0)
			\psline[border=0pt](1,0)(1,2)
			
			\psline[border=0pt]{-}(1.5,2)(1.5,0)
			\psline[border=0pt]{-}(1.5,0)(2.5,0)
			\psline[border=0pt]{-}(2.5,0)(2.5,2)

			\psline[border=0pt]{-}(3,2)(3,0)
			\psline[border=0pt]{-}(3,0)(4,0)
			\psline[border=0pt]{-}(4,0)(4,2)
			
			\psdots[dotsize=1.5pt](4.5,1)(5,1)(5.5,1)
			
			\psline[border=0pt]{-}(6,2)(6,0)
			\psline[border=0pt]{-}(6,0)(7,0)
			\psline[border=0pt]{-}(7,0)(7,2)
			
			\psline[border=0pt]{-}(7.5,2)(7.5,0)
			\psline[border=0pt]{-}(7.5,0)(8.5,0)
			\psline[border=0pt]{-}(8.5,0)(8.5,2)
			
			\rput(.5,-.5){{$1$}}
			\rput(2,-.5){{$2$}}
			\rput(3.5,-.5){{$3$}}
			\rput(6.5,-.5){{$U-1$}}
			\rput(8,-.5){{$U$}}
			
			\pscircle[fillstyle=solid,fillcolor=black, linewidth=1pt,linecolor=black](.3,.3){.2}
			\pscircle[fillstyle=solid,fillcolor=white, linewidth=1pt,linecolor=black](.7,.55){.2}
			\pscircle[fillstyle=solid,fillcolor=black, linewidth=1pt,linecolor=black](.3,.8){.2}
			
			\pscircle[fillstyle=solid,fillcolor=white, linewidth=1pt,linecolor=black](1.8,.3){.2}
			\pscircle[fillstyle=solid,fillcolor=white, linewidth=1pt,linecolor=black](2.2,.55){.2}
			\pscircle[fillstyle=solid,fillcolor=white, linewidth=1pt,linecolor=black](1.8,.8){.2}

			\pscircle[fillstyle=solid,fillcolor=white, linewidth=1pt,linecolor=black](3.3,.3){.2}
			\pscircle[fillstyle=solid,fillcolor=white, linewidth=1pt,linecolor=black](3.7,.55){.2}
			\pscircle[fillstyle=solid,fillcolor=black, linewidth=1pt,linecolor=black](3.3,.8){.2}

			\pscircle[fillstyle=solid,fillcolor=black, linewidth=1pt,linecolor=black](6.3,.3){.2}
			\pscircle[fillstyle=solid,fillcolor=black, linewidth=1pt,linecolor=black](6.7,.55){.2}
			\pscircle[fillstyle=solid,fillcolor=black, linewidth=1pt,linecolor=black](6.3,.8){.2}

			\pscircle[fillstyle=solid,fillcolor=black, linewidth=1pt,linecolor=black](7.8,.3){.2}
			\pscircle[fillstyle=solid,fillcolor=white, linewidth=1pt,linecolor=black](8.2,.55){.2}
			\pscircle[fillstyle=solid,fillcolor=white, linewidth=1pt,linecolor=black](7.8,.8){.2}
	\endpspicture

\noindent
{\footnotesize Figure 2:
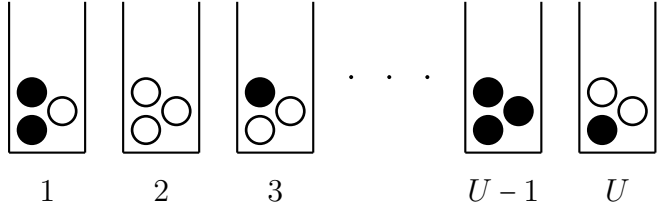
At $n=3$ each urn contains $3$ balls.}
\setcounter{figure}{2}

\vspace{0.3cm}
\noindent
	\end{center}

Let us denote  
by $B_n^u$ and $W_n^u$, respectively, the number of black and white balls in the urn $u$ at time $n$. 
Therefore $B_n^*:=\sum_uB_n^u$ (resp.~$W_n^*:=\sum_uW_n^u$) is
the total number of black (resp.~white) balls in the system at time $n$. 
Note that $B_n^*+W_n^*= Un$.

Define $\pi(i,j)=\frac{w_i}{w_i+w_{j}}$, and note that $\pi(i,j)$ is the probability of drawing
 a black ball out of an urn containing  $i+j$ balls among which $i$ are black, 
under the reinforcement mechanism determined by 
the weights $(w_i)_{i\in\N\cup\{0\}}$.

- The configuration at time $n+1$ is built from that at time $n$ according to the following transition
mechanism: 
Let $\xi_n^u$ be $U$ independent Bernoulli($p$) random variables, independent from the filtration up to time $n$.
For each $u\in \{1,\ldots, U\}$, 

\begin{itemize}
\item[{(a)}] On $\{\xi_n^u=1\}$, draw a ball out of all the $U$ urns combined.
The color of the drawn ball is black with probability
$\pi\left(B_n^*,W_n^*\right)$, and white with probability $\pi\left(W_n^*,B_n^*\right)$;
\item[{(b)}] On $\{\xi_n^u=0\}$, draw a ball out of urn $u$ alone,
so that the probability of black is $\pi\left(B_n^u,W_n^u\right)$, and of white is $\pi\left(W_n^u,B_n^u\right)$;
\item[{(c)}]
record the color of the drawn ball and return it to its urn of origin, 
and add another ball of the same color to urn $u$.
\end{itemize}
Doing the above procedure (a)--(c) simultaneously for each $u\in \{1,\ldots,U\}$,
 yields the configuration at time $n+1$. Call this process the {\em Interacting Urn Mechanism}.
\medskip

Note that when $p=0$ the urns are completely independent, while if $p=1$ the system behaves just like a single urn model in which $U$ balls are drawn with replacement and $U$ new balls that match in color are added in, at each step.

The special case $U=1$ and $w_i=1+i\Delta$ was firstly introduced by P\'{o}lya in 1930 to model the spreading of an epidemic (P\'{o}lya \cite{polya}). A number of generalizations have been introduced since, concerning the reinforcement sequence, the number of colors or the way of replacing balls into the urns. A good introduction about these models can be found in a survey by Pemantle \cite{pemantle_survey}.

An interesting condition to consider about the reinforcement weight sequence is the {\em Strong Reinforcement Hypothesis}:
\begin{equation}
\label{SRP}
\tag{SRH}
\sum_{i=0}^{\infty}\frac{1}{w_i}<\infty.
\end{equation}

Define event
$A:=\{\exists N<\infty$, such that starting from time $N$, all the balls drawn (added) have the same color$\}$.
The strongly reinforced mechanism is interesting due to the following result, proved for example in \cite{davis} or \cite{pemantle_survey} using Rubin's exponential embedding.

\begin{prop}
\label{p1}
If $U=1$ then
$$(\ref{SRP})\Leftrightarrow \mathbb{P}(A)=1. $$
\end{prop}

\textbf{Remark:} It is clear that the value of $p$ has no influence on 
the dynamics when $U=1$.

It is natural to ask whether the behavior of Proposition \ref{p1} persists under (\ref{SRP})
in the general case $U\geq2$ and $p\in [0,1]$. Though it seems to be a difficult question in the general case, we shall study in this paper the case of exponential reinforcement weight sequences, that is $w_i=\rho^i$ with $\rho>1$, which is much simpler.

The main theorem of this paper is the following.

\begin{theorem}
If the reinforcement weight sequence satisfies $\liminf_{i\to\infty}w_{i+1}/w_i>1$ then
\begin{itemize}
\item If $p\geq 1/2$ then all the urns eventually fixate on the same color ;
\item if $p<1/2$ then there exist a color $\mathfrak{c}$ such that eventually, some urns fixate on $\mathfrak{c}$ and the other urns draw a ball of color $\mathfrak{c}$ each time the drawing is done in all the urns combined and of the other color each time the drawing is done in the urn alone.
\end{itemize}
\end{theorem}

A more precise statement of these different cases will be given in section 3. In the major part of this paper we will limit ourselves to only two colors of balls (black and white) and to exponential reinforcement sequences. However we explain in the last section how the proofs extend easily to any number of colors and to any reinforcement weight sequence that satisfy $\liminf_{i\to\infty}w_{i+1}/w_i>1$.

It is worth noting that Interacting Urn Models is closely linked to multiple particles Reinforced Random Walks on star shaped graphs:

\begin{center}
	\pspicture(-1,-1)(1,1)
			\psline[linewidth=1pt]{*-*}(0,0)(1,0)
			\psline[linewidth=1pt]{*-*}(0,0)(-1,0)
			\psline[linewidth=1pt]{*-*}(0,0)(0,-1)
			\psline[linewidth=1pt]{*-*}(0,0)(0,1)
			\psline[linewidth=1pt]{*-*}(0,0)(.7,.7)
			\psline[linewidth=1pt]{*-*}(0,0)(.7,-.7)
			\psline[linewidth=1pt]{*-*}(0,0)(-.7,.7)
			\psline[linewidth=1pt]{*-*}(0,0)(-.7,-.7)
	\endpspicture

\noindent
{\footnotesize Figure 3:
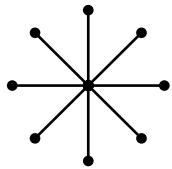
A star shaped graph with $8$ edges.}
\setcounter{figure}{3}

\vspace{0.3cm}
\noindent
	\end{center}
\noindent where the number of particles is equal to the number of urns and the number of edges is equal to the number of colors. The particles interact by sharing their memory just as the urns do.
In Reinforced Random Walks, (\ref{SRP}) also plays an important role since this is expected to be the condition for the walk getting eventually stuck on one edge. This assertion, known as the Sellke conjecture, has been proven in \cite{limictarres} by Limic and Tarr\`{e}s for all $(w_i)_{i\in\N\cup\{0\}}$ satisfying an additional technical condition, including nondecreasing $(w_i)_{i\in\N\cup\{0\}}$.
In another direction, if the reinforcement is affine, then it was proved by Enriquez and Sabot in \cite{sabot1} that some Reinforced Random Walks behave as an annealed Random Walk in Dirichlet Environment (see  \cite{sabot2}), Merkl and Rolles prove in \cite{merkl} the recurrence of such a walk in the graphs obtained from $\Z^2$ by replacing every edge by a sufficiently large, but fixed number of edges in series. In the case of polynomial reinforcements T\'{o}th proved in \cite{toth} that there are some limit theorems for weakly reinforced random walks on $\Z$.
\medskip

The rest of the paper is organised as follows. Section 2 introduces all the notations and definitions about Urn Paths and Interacting Urn Mechanism. Section 3 is devoted to the exponential reinforcement and splits into one subsection about heuristics that shows us how a phase transition appears for $p=1/2$, and three subsections devoted respectively to supercritical, critical and subcritical cases. Section 4 discusses two generalizations of the model: the first one shows that our results are still true if the number of colors is larger than $2$, and the second one extends them to a wider class of reinforcement weight sequences. The conclusion provides some open questions about the model.

\section{Definitions}

\subsection{\ups}

Denote by $\mathcal{U}=\{1,\cdots,U\}$ the set of $U$ urn labels, where $U\in\N$. Denote by $\mathfrak{C}=\{\mathfrak{b},\mathfrak{w}\}$ the set of two colors: black ($\mathfrak{b}$) and white ($\mathfrak{w}$).

\begin{defi}
We call a sequence $r=(r_n)_{n\in\N\cup\{0\}}$ on $\Z^{U}$ an {\em \up} if:
\begin{itemize}
	\item $r_0 =0_{\Z^{U}}$,
	\item $\forall u\in\mathcal{U}, \forall n\in\N$, $r_{n}(u)-r_{n-1}(u)\in\{0,1\}$.
\end{itemize}
And we denote by $\mathcal{R}$ the set of all {\ups}.
\end{defi}
%\up abreviates as UP (UP)
The coordinates $r_n(u)$ of $r_n$ should be interpreted as the numbers of black balls in urn $u$ at time $n$.

\begin{prop}
\label{conservn}
If $r$ is an \up, then
\begin{enumerate}[i)]
\item $\forall n\in\N, \forall u\in\mathcal{U},~ 0\leq r_n(u)\leq n;$
\item $s=(n-r_n)_{n\in\N\cup\{0\}}$ is also an {\up}. Its coordinates $s_n(u)=n-r_n(u)$ should be interpreted as the number of white balls in urn $u$ at time $n$.
\end{enumerate}
\end{prop}

\begin{proof}
The proof by induction is immediate from the definition of an \up.
\end{proof}

\begin{defi}
\label{defup}
Let $r\in\mathcal{R}$, $u\in\mathcal{U}$, and $n\in\N\cup\{0\}$ be a time unit. Let us denote by:
\begin{itemize}
	\item $\mathfrak{c}_n^u$ the color of the ball which is added to the urn $u$ at time $n\geq 1$:
$$\mathfrak{c}_n^u = \left\{
          \begin{array}{ll}
            \mathfrak{b} & \qquad \mathrm{if}\quad r_n(u)-r_{n-1}(u) = 1,\\
            \mathfrak{w} & \qquad \mathrm{if}\quad r_n(u)-r_{n-1}(u) = 0.\\
          \end{array}
        \right.$$
	\item $b_n^u$ the number of black balls in urn $u$ at time $n$, in symbols $b_n^u = r_n(u);$
	\item $w_n^u$ the number of white balls in urn $u$ at time $n$, in symbols $w_n^u = n-r_n(u);$
	\item $b_n^*$ the total number of black balls in the $U$ urns at time $n$: $b_n^* = \sum_{u\in\mathcal{U}}b_n^u = \sum_{u\in\mathcal{U}}r_n(u);$
	\item $w_n^*$ the total number of white balls in the $U$ urns at time $n$: $w_n^* = \sum_{u\in\mathcal{U}}w_n^u = Un-\sum_{u\in\mathcal{U}}r_n(u).$
	\item $\mathfrak{m}_n^u$ the majority color in urn $u$ at time $n$:
$$\mathfrak{m}_n^u = \left\{
          \begin{array}{ll}
            \mathfrak{b} & \qquad \mathrm{if}\quad b_n^u\geq w_n^u\\
            \mathfrak{w} & \qquad \mathrm{if}\quad b_n^u< w_n^u \\
          \end{array}
        \right..$$
\item $\mathfrak{m}_n^*$ the majority color in the $U$ urns at time $n$:
$$\mathfrak{m}_n^* = \left\{
          \begin{array}{ll}
            \mathfrak{b} & \qquad \mathrm{if}\quad b_n^*\geq w_n^*\\
            \mathfrak{w} & \qquad \mathrm{if}\quad b_n^*< w_n^* \\
          \end{array}
        \right..$$
\end{itemize}
In case of ties, we choose the black as the majority color, but this is out of convenience and it is not really important for the results. %an arbitrary choice
\end{defi}

Sometimes, it will be useful to consider a finer time indexing for which at each step a ball is added to only one urn.
We will reserve the variable $k$ and a superscript prime for this time indexing.
In particular, for all $n\in\N\cup\{0\}$ and $u\in\mathcal{U}$, we have $\mathfrak{c}_k^{\prime}=\mathfrak{c}_n^u$, $b_k^{\prime}=b_n^u$, $\mathfrak{m}_k^{\prime}=\mathfrak{m}_n^u$, $x_k^{\prime}=x_n^u$ and $X_k^{\prime}=X_n^u$ whenever $k\equiv k(n,u)=nU+u$ ($x_n^u$ and $X_n^u$ shall be define in the next section).

\begin{center}
	\pspicture(-1.5,-1)(10.8,6.5)

\rput(-1.5,5.5){{$n=0$}}
\rput(-1.5,3){{$n=1$}}
\rput(-1.5,.5){{$n=2$}}

\psline[border=1.5pt](0,6)(0,5)
			\psline[border=0pt](0,5)(1,5)
			\psline[border=0pt](1,5)(1,6)
			
			\psline[border=0pt]{-}(2.2,6)(2.2,5)
			\psline[border=0pt]{-}(2.2,5)(3.2,5)
			\psline[border=0pt]{-}(3.2,5)(3.2,6)

			\psline[border=0pt]{-}(4.4,6)(4.4,5)
			\psline[border=0pt]{-}(4.4,5)(5.4,5)
			\psline[border=0pt]{-}(5.4,5)(5.4,6)
			
			\psdots[dotsize=1.5pt](5.9,5.5)(6.4,5.5)(6.9,5.5)
			
			\psline[border=0pt]{-}(7.6,6)(7.6,5)
			\psline[border=0pt]{-}(7.6,5)(8.6,5)
			\psline[border=0pt]{-}(8.6,5)(8.6,6)
			
			\psline[border=0pt]{-}(9.8,6)(9.8,5)
			\psline[border=0pt]{-}(9.8,5)(10.8,5)
			\psline[border=0pt]{-}(10.8,5)(10.8,6)
			
			\rput(.5,4.5){{$k=1$}}
			\rput(2.7,4.5){{$k=2$}}
			\rput(4.9,4.5){{$k=3$}}
			\rput(8.1,4.5){{$k=U-1$}}
			\rput(10.3,4.5){{$k=U$}}

			\psline[border=0pt](0,3.5)(0,2.5)
			\psline[border=0pt](0,2.5)(1,2.5)
			\psline[border=0pt](1,2.5)(1,3.5)
			
			\psline[border=0pt]{-}(2.2,3.5)(2.2,2.5)
			\psline[border=0pt]{-}(2.2,2.5)(3.2,2.5)
			\psline[border=0pt]{-}(3.2,2.5)(3.2,3.5)

			\psline[border=0pt]{-}(4.4,3.5)(4.4,2.5)
			\psline[border=0pt]{-}(4.4,2.5)(5.4,2.5)
			\psline[border=0pt]{-}(5.4,2.5)(5.4,3.5)
			
			\psdots[dotsize=1.5pt](5.9,3)(6.4,3)(6.9,3)
			
			\psline[border=0pt]{-}(7.6,3.5)(7.6,2.5)
			\psline[border=0pt]{-}(7.6,2.5)(8.6,2.5)
			\psline[border=0pt]{-}(8.6,2.5)(8.6,3.5)
			
			\psline[border=0pt]{-}(9.8,3.5)(9.8,2.5)
			\psline[border=0pt]{-}(9.8,2.5)(10.8,2.5)
			\psline[border=0pt]{-}(10.8,2.5)(10.8,3.5)
			
			\rput(.5,2){{$k=U+1$}}
			\rput(2.7,2){{$k=U+2$}}
			\rput(4.9,2){{$k=U+3$}}
			\rput(8.1,2){{$k=2U-1$}}
			\rput(10.3,2){{$k=2U$}}
			
			\pscircle[fillstyle=solid,fillcolor=black, linewidth=1pt,linecolor=black](.3,2.8){.2}
						
			\pscircle[fillstyle=solid,fillcolor=white, linewidth=1pt,linecolor=black](2.5,2.8){.2}
			
			\pscircle[fillstyle=solid,fillcolor=white, linewidth=1pt,linecolor=black](4.7,2.8){.2}
			
			\pscircle[fillstyle=solid,fillcolor=black, linewidth=1pt,linecolor=black](7.9,2.8){.2}
						
			\pscircle[fillstyle=solid,fillcolor=black, linewidth=1pt,linecolor=black](10.1,2.8){.2}

			\psline[border=0pt](0,1)(0,0)
			\psline[border=0pt](0,0)(1,0)
			\psline[border=0pt](1,0)(1,1)
			
			\psline[border=0pt]{-}(2.2,1)(2.2,0)
			\psline[border=0pt]{-}(2.2,0)(3.2,0)
			\psline[border=0pt]{-}(3.2,0)(3.2,1)

			\psline[border=0pt]{-}(4.4,1)(4.4,0)
			\psline[border=0pt]{-}(4.4,0)(5.4,0)
			\psline[border=0pt]{-}(5.4,0)(5.4,1)
			
			\psdots[dotsize=1.5pt](5.9,.5)(6.4,.5)(6.9,.5)
			
			\psline[border=0pt]{-}(7.6,1)(7.6,0)
			\psline[border=0pt]{-}(7.6,0)(8.6,0)
			\psline[border=0pt]{-}(8.6,0)(8.6,1)
			
			\psline[border=0pt]{-}(9.8,1)(9.8,0)
			\psline[border=0pt]{-}(9.8,0)(10.8,0)
			\psline[border=0pt]{-}(10.8,0)(10.8,1)
			
			\rput(.5,-.5){{$k=2U+1$}}
			\rput(2.7,-.5){{$k=2U+2$}}
			\rput(4.9,-.5){{$k=2U+3$}}
			\rput(8.1,-.5){{$k=3U-1$}}
			\rput(10.3,-.5){{$k=3U$}}
			
			\pscircle[fillstyle=solid,fillcolor=black, linewidth=1pt,linecolor=black](.3,.3){.2}
			\pscircle[fillstyle=solid,fillcolor=white, linewidth=1pt,linecolor=black](.7,.55){.2}
						
			\pscircle[fillstyle=solid,fillcolor=white, linewidth=1pt,linecolor=black](2.5,.3){.2}
			\pscircle[fillstyle=solid,fillcolor=white, linewidth=1pt,linecolor=black](2.9,.55){.2}
			
			\pscircle[fillstyle=solid,fillcolor=white, linewidth=1pt,linecolor=black](4.7,.3){.2}
			\pscircle[fillstyle=solid,fillcolor=white, linewidth=1pt,linecolor=black](5.1,.55){.2}

			\pscircle[fillstyle=solid,fillcolor=black, linewidth=1pt,linecolor=black](7.9,.3){.2}
			\pscircle[fillstyle=solid,fillcolor=black, linewidth=1pt,linecolor=black](8.3,.55){.2}
			
			\pscircle[fillstyle=solid,fillcolor=black, linewidth=1pt,linecolor=black](10.1,.3){.2}
			\pscircle[fillstyle=solid,fillcolor=white, linewidth=1pt,linecolor=black](10.5,.55){.2}
			
	\endpspicture

\noindent
{\footnotesize Figure 4:
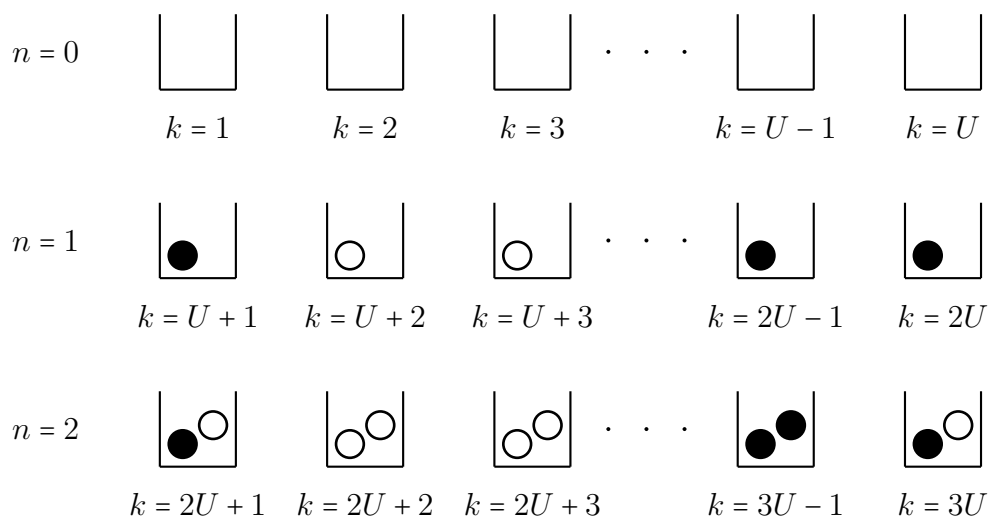
At time $k=nU+u$ one adds the $n$-th ball in the urn $u$.}
\setcounter{figure}{4}

\vspace{0.3cm}
\noindent
	\end{center}

\subsection{Interacting Urn Models}

Let $(\Omega,\mathcal{A},\mathbb{P})$ be a probability space on which we define:
\begin{enumerate}[i)]
\item a sequence $X=(X_n^u)_{n\in\N,u\in\mathcal{U}}$ of i.i.d random variables with distribution $p\delta_1+(1-p)\delta_{-1}$;
\item a sequence $V=(V_n^u)_{n\in\N,u\in\mathcal{U}}$ of i.i.d random variables with uniform distribution on $[0,1]$, independent of $X$.
\end{enumerate}
\begin{defi}
Denote by $\mathcal{F}_n$ the $\sigma$-field generated by the first $n$ steps, in symbols
$$\mathcal{F}_n=\sigma(X^u_m,V^u_m,~u\in\mathcal{U},~1\leq m\leq n).$$
\end{defi}

Fix the reinforcement weight sequence $w=(w_i)_{i\in\N\cup\{0\}}\in\mathbb{R}_+^{\N\cup\{0\}}$ and recall that $\pi(i,j)=\frac{w_i}{w_i+w_j}$ for $i,j\in\N\cup\{0\}$.
Set $\chi=\{-1,+1\}^{\mathcal{U}\times\mathbb{N}}$ and call $x=(x_n^u)_{n\in\N,u\in\mathcal{U}}\in\chi$ an environment.
We next define a process $(R_n)_{n\in\N}$ on $\Z^U$ such that $R_0 = 0_{\Z^U}$ and $\forall n\in\N, \forall u\in\mathcal{U}$:

\begin{itemize} 
	\item if $x_n^u = 1$, then:
\begin{equation}
\tag{$\pentagram$}
\label{penta}
R_n(u)-R_{n-1}(u)=
\left\{
	\begin{array}{ll}
            1 & \text{if } V_n^u<\pi(B_{n-1}^*,W_{n-1}^*), \\
            0 & \text{if } V_n^u\geq\pi(B_{n-1}^*,W_{n-1}^*); \\
          \end{array}
\right.
\end{equation}

\item if $x_n^u = -1$, then:
\begin{equation}
\tag{$\pentagram\pentagram$}
\label{2penta}
R_n(u)-R_{n-1}(u)=
\left\{
	\begin{array}{ll}
            1 & \text{if } V_n^u<\pi(B_{n-1}^u,W_{n-1}^u), \\
            0 & \text{if } V_n^u\geq\pi(B_{n-1}^u,W_{n-1}^u). \\
          \end{array}
\right.
\end{equation}
\end{itemize}

Note that $R$ is almost surely an {\up}. Hence we can define the quantities $B_n^u$, $W_n^u$, $B_n^*$ and $W_n^*$ for all $n\geq 0$ and $u\in\mathcal{U}$ as in Definition \ref{defup}. This justifies the use of these notations in (\ref{penta}) and (\ref{2penta}).

We call the law of this process the {\em Interacting Urn Mechanism in the environment $x$} and denote it by $\text{IUM}^x$. We denote by $P^x$ the probability it induces on the set of all \ups. 

Now fix the interaction parameter $p\in[0,1]$. If the environment is sampled as a sequence $(X_n^u)_{n\in\N,u\in\mathcal{U}}$ of i.i.d random variables such that:
$$\mathbb{P}\left[X_n^u=1\right]=1-\mathbb{P}\left[X_n^u=-1\right]=p, ~~\forall n\in\N, \forall u\in\mathcal{U},$$
then this process is called the {\em{Interacting Urn Mechanism}}. We denote its law by $\RUPa$, and denote by $P$ the probability it induces on the set of all \ups.

\begin{defi}
\label{defi_G}
Let $\mathcal{G}_1$ be the $\sigma$-field generated by what happens when the balls are drawn out of {\drawall}, in symbols
$$\mathcal{G}_1=\sigma\left((X_n^u)_{n\in\N,u\in\mathcal{U}},(V_n^u\mathbb{1}_{\{X_n^u=1\}})_{n\in\N,u\in\mathcal{U}}\right),$$
and $\mathcal{G}_{-1}$ the $\sigma$-field generated by what happens when the balls are drawn out of one urn at a time, in symbols
$$\mathcal{G}_{-1}=\sigma\left((X_n^u)_{n\in\N,u\in\mathcal{U}},(V_n^u\mathbb{1}_{\{X_n^u=-1\}})_{n\in\N,u\in\mathcal{U}}\right).$$
\end{defi}

\section{The exponential reinforcement}

From now on we shall focus on the case where $w_i=\rho^i$ for some $\rho>1$. That is to say the exponential {\RUPa} that satisfies (\ref{SRP}).

The important thing to notice about exponential reinforcement is that the transition law depends only on the difference $B-W$:

$$\pi(B,W)=\frac{\rho^{B}}{\rho^{B}+\rho^{W}}=\frac{1}{1+\rho^{W-B}},$$

and

$$\pi(W,B)=\frac{\rho^{W}}{\rho^{B}+\rho^{W}}=\frac{1}{1+\rho^{B-W}}.$$

\subsection{Heuristics}

Let us consider initially the case $U=2$ and assume that the probability of drawing a black ball out of urn $u$ at time $n$ converges to a limit as $n$ goes to the infinity (which we shall prove in the next subsections). Denote this limit by $\alpha_u$:
$$\alpha_u:=\lim_{n\rightarrow\infty}\mathbb{P}\left(\mathfrak{b} \text{ is drawn out of urn }u\text{ at time }n\right), u=1\text{ or }2.$$

Then by a Law of Large Numbers we get:

$$\lim_{n\rightarrow\infty}\frac{B_n^u}{n}=\alpha_u, ~~u=1\text{ or }2.$$

Thus,

\begin{eqnarray*}
\alpha_u & = & \lim_{n\rightarrow\infty}\mathbb{P}\left(\mathfrak{b} \text{ is drawn out of } u\text{ at time }n\right)=\mathbb{P}\left[\mathfrak{c}_n^u=\mathfrak{b}\right]\\
&=&\lim_{n\rightarrow\infty}\left[\mathbb{P}\left[X_n^u=-1\right]\mathbb{P}\left[\left.\mathfrak{c}_n^u=\mathfrak{b}\right|X_n^u=-1\right]+\mathbb{P}\left[X_n^u=1\right]\mathbb{P}\left[\left.\mathfrak{c}_n^u=\mathfrak{b}\right|X_n^u=1\right]\right]\\
&=&\lim_{n\rightarrow\infty} (1-p)\frac{\rho^{B_n^u}}{\rho^{B_n^u}+\rho^{n-B_n^u}}
+p\frac{\rho^{B_n^1+B_n^2}}{\rho^{B_n^1+B_n^2}+\rho^{2n-(B_n^1+B_n^2)}}\\
&=&\lim_{n\rightarrow\infty} (1-p)\frac{1}{1+\rho^{n-2B_n^u}}
+p\frac{1}{1+\rho^{2n-2(B_n^1+B_n^2)}}\\
&=&\lim_{n\rightarrow\infty} (1-p)\frac{1}{1+\rho^{n(1-2\alpha_u)}}+p\frac{1}{1+\rho^{2n(1-\alpha_1-\alpha_2)}}\\
&=&(1-p)\mathbb{1}_{\left\{\alpha_u>\frac{1}{2}\right\}}+p\mathbb{1}_{\left\{\alpha_1+\alpha_2>1\right\}}.
\end{eqnarray*}

This yields the following system of two equations with two variables:

$$
\left\{
\begin{array}{l}
\alpha_1=(1-p)\mathbb{1}_{\left\{\alpha_1>\frac{1}{2}\right\}}+p\,\mathbb{1}_{\left\{\alpha_1+\alpha_2>1\right\}},\\
\alpha_2=(1-p)\mathbb{1}_{\left\{\alpha_2>\frac{1}{2}\right\}}+p\,\mathbb{1}_{\left\{\alpha_1+\alpha_2>1\right\}}.
\end{array}
\right.
$$

There are six possible cases to consider:

\begin{itemize}
\item If $\alpha_1<1/2$ and $\alpha_2<1/2$ (and therefore $\alpha_1+\alpha_2<1$), then the system above gives: $\alpha_1=\alpha_2=0$.
\item If $\alpha_1>1/2$ and $\alpha_2>1/2$ (and therefore $\alpha_1+\alpha_2>1$), then the system above gives: $\alpha_1=\alpha_2=(1-p)+p=1$.
\item If $\alpha_1<1/2$, $\alpha_2>1/2$ and $\alpha_1+\alpha_2>1$, then the system above gives: $\alpha_1=p$ and $\alpha_2=1$. This case, as well as the next three cases, is clearly possible only if $p<1/2$.
\item If $\alpha_1>1/2$, $\alpha_2<1/2$ and $\alpha_1+\alpha_2>1$, then the system above gives: $\alpha_1=1$ and $\alpha_2=p$.
\item If $\alpha_1<1/2$, $\alpha_2>1/2$ and $\alpha_1+\alpha_2<1$, then the system above gives: $\alpha_1=0$ and $\alpha_2=1-p$.
\item If $\alpha_1>1/2$, $\alpha_2<1/2$ and $\alpha_1+\alpha_2<1$, then the system above gives: $\alpha_1=1-p$ and $\alpha_2=0$.
\end{itemize}

This heuristic indicates that a phase transition occurs at $p=1/2$. We obtain the following configuration graph:

	\begin{center}
	\pspicture(-1,-1)(5,5)
			\psline[border=0pt]{->}(0,0)(0,4.5)
			\psline[border=0pt]{->}(0,0)(4.5,0)
			\psline[border=0pt]{-}(-.1,4)(.1,4)
			\psline[border=0pt]{-}(4,-.1)(4,.1)
			\psline[border=0pt,linestyle=dashed](0,1)(4.2,1)
			\psline[border=0pt,linestyle=dashed](1,0)(1,4.2)			
			\psdots[linecolor=black, dotsize=5pt](0,0)(4,4)
			\psdots[linecolor=gray, dotsize=5pt](3,0)(4,1)(0,3)(1,4)
			\rput(-.5,4){1}
			\rput(4,-.5){1}
			\rput(-.3,-.3){0}
			\rput(1,-.5){$p$}
			\rput(3,-.5){$1-p$}
			\rput(-.5,1){$p$}
			\rput(-.7,3){$1-p$}
			\rput(-.2,4.6){{$\alpha_2$}}
			\rput(4.6,-.2){{$\alpha_1$}}			
	\endpspicture
	
\noindent
{\footnotesize Figure 5:
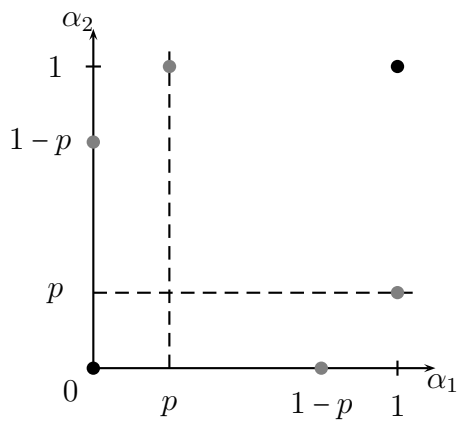
Graph of possible configurations when $U=2$.\\The two black points exists for any $p\in [0,1]$, the four gray points exists only for $p<1/2$.}
\setcounter{figure}{5}

\vspace{0.3cm}
\noindent
	\end{center}

\textbf{Remark:} Note that the solutions are symmetric with respect to the line $\alpha_1=\alpha_2$. This comes from the fact that inverting the labels of the urns gives a process with the same law. Similarly, inverting the two colors is manifested on the graph by a symmetry with respect to the point $\alpha_1=\alpha_2=1/2$. (For example, $(1-p,0)$ is symmetric to $(p,1)$)

Taking these two transformations into account show us that the two points $(0,0)$ and $(1,1)$ corespond in fact to two similar asymptotic behaviors, and the four points $(0,1-p)$, $(p,1)$, $(1-p,0)$ and $(1,p)$ are also similar. In particular, similar points occur with the same probability.

When $p>1/2$ we have only two possibilities (black points) and each of them has probability $1/2$.

When $p<1/2$ we have six possibilities. An interesting question is: what is the respective probability of each of the above points?
\\

One can object here that we forgot to discuss the case $\alpha_1 =1/2$ or $\alpha_2 =1/2$. Taking this possibility into account yields us to many other points on the configuration graph, (in fact, all the points on the lines $\alpha_1=1/2$ and $\alpha_2=1/2$). But since the evolution of the process depends only on the
differences $B_n^u-W_n^u$, we can conjecture that the cases when $\alpha_1$ or $\alpha_2$ equals to $1/2$ will be ``unstable solutions''. We shall
study this in more detail in the next subsections.

\subsection{The supercritical case}

In this subsection we consider the case when $p>1/2$.
We first state a simple but useful lemma that is a consequence of the law of large numbers and is left to the reader.

\begin{lemme}
\label{lemme1}
Let $Y_1, Y_2, Y_3 \dots$ be a sequence of i.i.d. random variables such that $\mathbb{P}[Y_1=1]=1-\mathbb{P}[Y_1=-1]=p$, then
$\forall \alpha\in(0,2p-1)$, we have:
$$\mathbb{P}\left[Y_1+\cdots +Y_k \geq \alpha k, \forall k\geq 0\right]>0.$$
\end{lemme}

Set $\theta_p$ to be the above probability. (Note that $\theta_p$ depends on $\alpha$.)

\begin{theorem}
\label{psur}
Assume that $p>1/2$.
$$\mathbb{P}\left[\exists n_0\in\N, \exists \mathfrak{d}\in\mathfrak{C} \text{ such that } \forall n\geq n_0, \forall u\in\mathcal{U} \text{ we have } \mathfrak{c}_n^u = \mathfrak{d}\right]=1$$
In words, after some finite time $n_0$ all the balls drawn have the same color $\mathfrak{d}$.
\end{theorem}

The theorem is a consequence of the following lemma.

%Recall that $\mathcal{A}(n,u)=\mathcal{A}^{\prime}(nU+u)$ is the event \textit{``The draw for urn $u$ at time $n$ is out of \drawall''} and $\mathcal{A}^c(n,u)=\mathcal{A}^{\prime c}(nU+u)$ the event \textit{``The draw for urn $u$ at time $n$ is out of the urn $u$ alone''}.

Recall from Definition \ref{defi_G} that $\mathcal{G}_1$ is the $\sigma$-field generated by the information of when the balls will be drawn out of {\drawall} ($(X_n^u)_{n\in\mathbb{N}, u_in\mathcal{U}}$) and of the parameters $(V_n^u\,\mathbb{1}_{\{(X_n^u)=1\}})_{n\in\mathbb{N}, u\in\mathcal{U}}$ that decide what happens when the balls are drawn out of {\drawall}.

\begin{lemme}\label{lemmesur1}
Assume that $p> 1/2$. There exists a constant $\eta_1>0$ such that for all $n_1\in\N$,
$$\mathbb{P}\left[\mathbb{P}\left[\bigcap_{n> n_1}\bigcap_{u\in\mathcal{U}}\{X_n^u=-1\}\cup \{\mathfrak{c}_n^u=\mathfrak{m}_{n_1}^*\}\Big|\mathcal{G}_1\right]=1 \Big| \mathcal{F}_{n_1}\right]\geq\eta_1~~a.s.$$
This means that with a probability uniformly bounded away from $0$, the information $\mathcal{G}_1$ is sufficient to be sure that after time $n_1$, if a ball is drawn out of {\drawall} its color is the majority color at time $n_1$.
\end{lemme}

\textbf{Remark: } Note that the event 
$$\bigcap_{n> n_1}\bigcap_{u\in\mathcal{U}}\{X_n^u=-1\}\cup \{\mathfrak{c}_n^u=\mathfrak{m}_{n_1}^*\}$$
is not $\mathcal{G}_1$ measurable since for any $n$ and $u$ the event $\{\mathfrak{c}_n^u=\mathfrak{m}_{n_1}^*\}$ is determined by comparison of the variable $V_n^u$ with the variable $\pi(B_{n-1}^*,W_{n-1}^*)$ which is not $\mathcal{G}_1$-measurable (c.f. the definition of the {{\RUPa}}.)
\begin{proof}
Without loss of generality, in this proof we will suppose that the majority color at time $n_1$ is black, in symbols, $\mathfrak{m}_{n_1}^*=\mathfrak{b}$.

For all $n\in\N$ and $u\in \mathcal{U}$ recall that $X_k^{\prime}:=X_n^u$ when $k=nU+u$ and set
$$Y_{k}=\mathbb{1}_{\{X_k^{\prime}=1\}}-\mathbb{1}_{\{X_k^{\prime}=-1\}}.$$

Set $k_1=(n_1+1)U$. The sequence $Y_{k_1+1}, Y_{k_1+2}, Y_{k_1+3},\cdots$ satisfies the conditions of Lemma \ref{lemme1} and is independent of $\mathcal{F}_{n_1}$. So for all $\alpha\in(0,2p-1)$ we have:
$$\mathbb{P}\left[A|\mathcal{F}_{n_1}\right]=\mathbb{P}\left[A\right]=\theta_p>0,$$
where $A= \{Y_{k_1+1}+ \cdots+Y_k\geq\alpha (k-k_1),\forall k\geq k_1\}$. Let us fix $\alpha=(2p-1)/2$ (any $\alpha\in (0,2p-1)$ would suffice) and set
\begin{equation}
\label{1etoile}
\gamma_n=\alpha U(n-n_1-1)
\end{equation}
and denote by $C^{\prime}_k=C^u_n$ the event
\begin{equation}
\label{2etoile}
C^{\prime}_k=C^u_n=\left\{V_n^u\,\mathbb{1}_{\{X_n^u=1\}}\leq\frac{1}{1+\rho^{-\gamma_{n}}}\right\}.
\end{equation}
Now set
$$C = A\cap\bigcap_{k=k_1+1}^{\infty}C^{\prime}_k.$$

The proof will proceed in two steps:
\begin{enumerate}[(i)]
\item \label{point1} We have the following event inclusion:
$$C\subset \left\{\mathbb{P}\left[\bigcap_{n> n_1}\bigcap_{u\in\mathcal{U}}\{X_n^u=-1\}\cup \{\mathfrak{c}_n^u=\mathfrak{m}_{n_1}^*\}\Big|\mathcal{G}_1\right]=1\right\}.$$
\item \label{point2} The probability of $C$ knowing $\mathcal{F}_{n_1}$ is uniformly bounded away from $0$:
$$\mathbb{P}\left[C \big| \mathcal{F}_{n_1}\right]\geq\eta_1~~ a.s.$$
\end{enumerate}

It is easy to see that the lemma follows from those two points. Namely,
$$0<\eta_1\leq\mathbb{P}\left[C|\mathcal{F}_{n_1}\right]\leq\mathbb{P}\left[\mathbb{P}\left[\bigcap_{n> n_1}\bigcap_{u\in\mathcal{U}}\{X_n^u=-1\}\cup \{\mathfrak{c}_n^u=\mathfrak{m}_{n_1}^*\}\Big|\mathcal{G}_1\right]=1 \Big| \mathcal{F}_{n_1}\right]~a.s.$$

Let us prove (\ref{point1}).
First, we prove by induction that:

$$C\subset \bigcap_{n> n_1}\bigcap_{u\in\mathcal{U}}\{X_n^u=-1\}\cup \{\mathfrak{c}_n^u=\mathfrak{m}_{n_1}^*\}.$$

Suppose that we know
\begin{equation}
\label{4etoile}
C\subset \bigcap_{n=n_1+1}^N\bigcap_{u\in\mathcal{U}}\{X_n^u=-1\}\cup \{\mathfrak{c}_n^u=\mathfrak{m}_{n_1}^*\},
\end{equation}

which means that on $C$, each time a ball is drawn out of {\drawall} between time $n_1+1$ and $N$, then it is black. So (\ref{4etoile}) gives us:

\begin{eqnarray*}
B_{N}^*-W_{N}^* &\geq& (B_{n_1}^*-W_{n_1}^*)+\sum_{n=n_1+1}^N\sum_{u\in\mathcal{U}}\mathbb{1}_{\{X_n^u=1\}} - \sum_{n=n_1+1}^N\sum_{u\in\mathcal{U}}\mathbb{1}_{\{X_n^u=-1\}}\\
&\geq&\sum_{n=n_1+1}^N\sum_{u\in\mathcal{U}}\mathbb{1}_{\{X_n^u=1\}} - \sum_{n=n_1+1}^N\sum_{u\in\mathcal{U}}\mathbb{1}_{\{X_n^u=-1\}}\\
&\geq&\sum_{k=k_1+1}^{U(N+1)}\left(\mathbb{1}_{\{X_k^\prime=1\}}-\mathbb{1}_{\{X_k^\prime=-1\}}\right)\\
&\geq&\sum_{k=k_1+1}^{U(N+1)}Y_k.
\end{eqnarray*}
And since $C\subset A$ by the definition of $A$ for $k=U(N+1)$:

\begin{equation}
\label{3etoile}
B_{N}^*-W_{N}^*\geq \alpha(k-k_1)=\alpha(U(N+1)-U(n_1+1))=\alpha U (N-n_1) = \gamma_{N+1},
\end{equation}
recalling (\ref{1etoile}). So for each $u\in\mathcal{U}$, using (\ref{2etoile}) with $n=N+1$,
we deduce that on $C$, $V_{N+1}^u \leq\frac{1}{1+\rho^{-\gamma_{N+1}}} \leq\frac{1}{1+\rho^{-(B_N^*-W_N^*)}}$ due to (\ref{3etoile}). Due to the construction of the {{\RUPa}} (\ref{penta}) the last observation means that if $X_n^u=1$ then the ball added to $u$ at time $N+1$ is black.

It follows that

$$C\subset \{X_{N+1}^u=-1\}\cup \{\mathfrak{c}_{N+1}^u=\mathfrak{m}_{n_1}^*=\mathfrak{b}\},~~~\forall u\in\mathcal{U},$$

and therefore by induction that

$$C\subset \bigcap_{n> n_1}^{N+1}\bigcap_{u\in\mathcal{U}}\{X_n^u=-1\}\cup \{\mathfrak{c}_n^u=\mathfrak{m}_{n_1}^*\}.$$

It remains to note that $C$ is $\mathcal{G}_1$-measurable, to conclude:

$$C = \left\{\mathbb{P}\left[C\big|\mathcal{G}_1\right]=1\right\} \subset \left\{\mathbb{P}\left[\bigcap_{n\geq n_1}\bigcap_{u\in\mathcal{U}}\{X_n^u=-1\}\cup \{\mathfrak{c}_n^u=\mathfrak{m}_{n_1}^*\}\Big|\mathcal{G}_1\right]=1\right\},$$
which shows (\ref{point1}).
Let us now prove (\ref{point2}). We have

\begin{eqnarray*}
\mathbb{P}\left[C \big| \mathcal{F}_{n_1}\right]&=&\mathbb{P}\left[A\cap\bigcap_{k=k_1+1}^{\infty}C^{\prime}_k \big| \mathcal{F}_{n_1}\right]\\
&=&\mathbb{E}\left[\mathbb{1}_A\times\prod_{k=k_1+1}^{\infty}\mathbb{1}_{C^{\prime}_k} \big| \mathcal{F}_{n_1}\right]\\
&=&\mathbb{E}\left[\mathbb{1}_A\times\mathbb{E}\left[\prod_{k=k_1+1}^{\infty}\mathbb{1}_{C^{\prime}_k}\Big| \mathcal{F}_{n_1};A \right]\big| \mathcal{F}_{n_1}\right].
\end{eqnarray*}

We next give a lower bound for the inside expectation. Fix $K> k_1$ and set $N=\left\lfloor\frac{K-1}{U}\right\rfloor$ and $\mathcal{F}_{n_1}^K=\sigma(\mathcal{F}_{n_1};A,C^{\prime}_{k_1+1},...,C^{\prime}_{K-1})$. Then, on $A$

\begin{eqnarray*}
\mathbb{E}\left[\prod_{k=k_1+1}^{K}\mathbb{1}_{C^{\prime}_k}\big| \mathcal{F}_{n_1};A \right]&=&\mathbb{E}\left[\mathbb{E}\left[\mathbb{1}_{C^{\prime}_K}\Big|\mathcal{F}_{n_1}^K\right]\times\prod_{k=k_1+1}^{K-1}\mathbb{1}_{C^{\prime}_k}\bigg| \mathcal{F}_{n_1};A \right]\\
&=&\mathbb{E}\left[\mathbb{E}\left[\mathbb{1}_{\{X_K^{\prime}=-1\}\cup\{V_K^{\prime} \leq\frac{1}{1+\rho^{-\gamma_N}}\}}\Big|\mathcal{F}_{n_1}^K\right]\times\prod_{k=k_1+1}^{K-1}\mathbb{1}_{C^{\prime}_k}\bigg| \mathcal{F}_{n_1};A \right],\\
\end{eqnarray*}
where due to the fact that $V_K^{\prime}$ is independent of $\mathcal{F}_{n_1}$, $A$, $C^{\prime}_{k_1+1}$,..., $C^{\prime}_{K-1}$ and $X_K^{\prime}$
\begin{eqnarray*}
\mathbb{E}\left[\mathbb{1}_{\{X_K^{\prime}=-1\}\cup\{V_K^{\prime} \leq\frac{1}{1+\rho^{-\gamma_N}}\}}\Big|\mathcal{F}_{n_1}^K\right]&=&\mathbb{E}\left[\mathbb{1}_{\{X_K^{\prime}=-1\}}+\mathbb{1}_{\{X_K^{\prime}=1\}\cap\{V_K^{\prime} \leq\frac{1}{1+\rho^{-\gamma_N}}\}}\Big|\mathcal{F}_{n_1}^K\right]\\
&=&\mathbb{E}\left[\mathbb{E}\left[\mathbb{1}_{\{X_K^{\prime}=-1\}}+\mathbb{1}_{\{X_K^{\prime}=1\}\cap\{V_K^{\prime} \leq\frac{1}{1+\rho^{-\gamma_N}}\}}\big|\mathcal{F}_{n_1}^K;X_K^{\prime}\right]\Big|\mathcal{F}_{n_1}^K\right]\\
&=&\mathbb{E}\left[\mathbb{1}_{\{X_K^{\prime}=-1\}}+\mathbb{1}_{\{X_K^{\prime}=1\}}\mathbb{E}\left[\mathbb{1}_{\{V_K^{\prime} \leq\frac{1}{1+\rho^{-\gamma_N}}\}}\big|\mathcal{F}_{n_1}^K;X_K^{\prime}\right]\Big|\mathcal{F}_{n_1}^K\right]\\
&=&\mathbb{E}\left[\mathbb{1}_{\{X_K^{\prime}=-1\}}+\mathbb{1}_{\{X_K^{\prime}=1\}}\frac{1}{1+\rho^{-\gamma_N}}\Big|\mathcal{F}_{n_1}^K\right]\\
&\geq&\mathbb{E}\left[(\mathbb{1}_{\{X_K^{\prime}=-1\}}+\mathbb{1}_{\{X_K^{\prime}=1\}})\frac{1}{1+\rho^{-\gamma_N}}\Big|\mathcal{F}_{n_1}^K\right]\\
&\geq&\mathbb{E}\left[\frac{1}{1+\rho^{-\gamma_N}}\Big|\mathcal{F}_{n_1}^K\right]\geq\frac{1}{1+\rho^{-\gamma_N}}.
\end{eqnarray*}
Therefore

\begin{eqnarray*}
\mathbb{E}\left[\prod_{k=k_1+1}^{K}\mathbb{1}_{C^{\prime}_k}\big| \mathcal{F}_{n_1};A \right]&\geq&\mathbb{E}\left[\frac{1}{1+\rho^{-\gamma_N}}\times\prod_{k=k_1+1}^{K-1}\mathbb{1}_{C^{\prime}_k}\bigg| \mathcal{F}_{n_1};A \right]\\
&\geq&\frac{1}{1+\rho^{-\gamma_N}}\mathbb{E}\left[\prod_{k=k_1+1}^{K-1}\mathbb{1}_{C^{\prime}_k}\bigg| \mathcal{F}_{n_1};A \right],
\end{eqnarray*}

and by induction we can conclude that

$$
\mathbb{E}\left[\prod_{k=k_1+1}^{K}\mathbb{1}_{C^{\prime}_k}\big| \mathcal{F}_{n_1};A \right]\geq\prod_{k=k_1+1}^{K-1}\frac{1}{1+\rho^{-\gamma_{n(k)}}}\geq\prod_{k=k_1+1}^{\infty}\frac{1}{1+\rho^{-\gamma_{n(k)}}},
$$

with $n(k) = \lfloor\frac{k-1}{U}\rfloor$. Since the last inequality is true for any $K> k_1$ we can let $K$ go to infinity:

$$
\mathbb{E}\left[\prod_{k=k_1+1}^{\infty}\mathbb{1}_{C^{\prime}_k}\big| \mathcal{F}_{n_1};A \right]\geq\prod_{k=k_1+1}^{\infty}\frac{1}{1+\rho^{-\gamma_{n(k)}}}=\prod_{n=n_1+1}^{\infty}\left(\frac{1}{1+\rho^{-\gamma_n}}\right)^U=\prod_{n=0}^{\infty}\left(\frac{1}{1+\rho^{-\alpha U n}}\right)^U.
$$

Hence we have

\begin{eqnarray*}
\mathbb{P}\left[C \big| \mathcal{F}_{n_1}\right]&=&\mathbb{E}\left[\mathbb{1}_A\times\mathbb{E}\left[\prod_{k=k_1+1}^{\infty}\mathbb{1}_{C^{\prime}_k}\Big| \mathcal{F}_{n_1};A \right]\Big| \mathcal{F}_{n_1}\right]\\
&\geq&\mathbb{E}\left[\mathbb{1}_A\times\prod_{n=0}^{\infty}\left(\frac{1}{1+\rho^{-\alpha U n}}\right)^U\Big| \mathcal{F}_{n_1}\right]\\
&\geq&\mathbb{E}\left[\mathbb{1}_A\big| \mathcal{F}_{n_1}\right]\times\prod_{n=0}^{\infty}\left(\frac{1}{1+\rho^{-\alpha U n}}\right)^U\\
&\geq&\theta_p\times\prod_{n=0}^{\infty}\left(\frac{1}{1+\rho^{-\alpha U n}}\right)^U.\\
\end{eqnarray*}

Therefore, setting
$$\eta_1=\theta_p\times\prod_{n=0}^{\infty}\left(\frac{1}{1+\rho^{-\alpha U n}}\right)^U>0,$$
we obtain (\ref{point2}) which gives the lemma.
\end{proof}

\begin{proof} (of Theorem \ref{psur})
Denote by $F_{n_1}$ the event
$$F_{n_1}=\left\{\mathbb{P}\left[\bigcap_{n\geq n_1}\bigcap_{u\in\mathcal{U}}\{X_n^u=-1\}\cup \{\mathfrak{c}_n^u=\mathfrak{m}_{n_1}^*\}\Big|\mathcal{G}_1\right]=1\right\}.$$
Note that $(F_{n_1})_{n_1\in\mathbb{N}}$ is an increasing sequence of events and let $F$ be the event
$$F=\lim_{n_1\rightarrow\infty}F_{n_1}=\bigcup_{n_1\in\mathbb{N}}F_{n_1}$$
then recall the L\'{e}vy 0-1 law (see for example \cite{durrett}):
$$\mathbb{1}_F=\lim_{n_1\rightarrow\infty}\mathbb{E}\left[\mathbb{1}_F|\mathcal{F}_{n_1}\right]=\lim_{n_1\rightarrow\infty}\mathbb{P}\left[F|\mathcal{F}_{n_1}\right]\geq\lim_{n_1\rightarrow\infty}\mathbb{P}\left[F_{n_1}|\mathcal{F}_{n_1}\right]\geq\eta_1>0\text{ a.s.}$$
by Lemma \ref{lemmesur1}. Therefore $\mathbb{P}(F)=1$. Denote by $n_2$ the random finite number defined by:
$$n_2=\min{\left\{n_1\in\mathbb{N} : \mathbb{P}\left[\bigcap_{n\geq n_1}\bigcap_{u\in\mathcal{U}}\{X_n^u=-1\}\cup \{\mathfrak{c}_n^u=\mathfrak{m}_{n_1}^*\}\Big|\mathcal{G}_1\right]=1\right\}}.$$
It is important to note here that $n_2$ is $\mathcal{G}_1$-measurable.

As in the proof of Lemma \ref{lemmesur1}, suppose without loss of generality that the above color of fixation for the $U$ urns combined is black ($\mathfrak{b}$).

Therefore, for a given $\beta\in(1/2,p)$, there is another finite time $n_3\geq n_2$ such that,
for all $n\geq n_3$ and $u\in\mathcal{U}$,
$$\frac{B_n^u}{n}\geq \beta,$$

and $n_3$ is $\mathcal{G}_1$-measurable. For example choose

$$n_3 = \min\left\{n\geq n_2 : \forall l>n,~\frac{\sum_{m=n_2}^l\mathbb{1}_{\{X^u_m=1\}}}{l}\geq \beta\right\}.$$
This minimum exists because
$$\lim_{l\to\infty}\frac{\sum_{m=n_2}^l\mathbb{1}_{\{X^u_m=1\}}}{l}=p> \beta$$
and it satisfies the expected property:
$$\forall n\geq n_3,~\frac{B_n^u}{n}=\frac{\sum_{m=1}^n\mathbb{1}_{\{c_n^u=\mathfrak{b}\}}}{n}\geq\frac{\sum_{m=n_2}^n\mathbb{1}_{\{c_n^u=\mathfrak{b}\}}}{n}\geq\frac{\sum_{m=n_2}^n\mathbb{1}_{\{X_n^u=1\}}}{n}\geq\beta.$$
Now let us bound the expected number of white balls in the urn $u$ as follows: firstly
$$\sum_{n=1}^\infty\mathbb{1}_{\{\mathfrak{c}_n^u={\mathfrak{w}}\}}= \sum_{n=1}^{n_3-1}\mathbb{1}_{\{\mathfrak{c}_n^u={\mathfrak{w}}\}}   + \sum_{n=n_3}^\infty\mathbb{1}_{\{\mathfrak{c}_n^u={\mathfrak{w}}\}} \leq n_3-1  + \sum_{n=n_3}^\infty\mathbb{1}_{\{\mathfrak{c}_n^u={\mathfrak{w}}\}},$$

and secondly

\begin{eqnarray*}
\mathbb{E}\left[\sum_{n=n_3}^\infty\mathbb{1}_{\{\mathfrak{c}_n^u={\mathfrak{w}}\}}\Big|\mathcal{G}_1\right]&=&\mathbb{E}\left[\sum_{n=n_3}^\infty\mathbb{1}_{\{\mathfrak{c}_n^u={\mathfrak{w}}\}\cap\{X_n^u=-1\}}\Big|\mathcal{G}_1\right]\\
&=&\mathbb{E}\left[\sum_{n=n_3}^\infty\mathbb{1}_{\{V_n^u>\frac{1}{1+\rho^{W^u_n-{B^u_n}}}\}\cap\{X_n^u=-1\}}\Big|\mathcal{G}_1\right]\\
&=&\sum_{n=n_3}^\infty\mathbb{E}\left[\mathbb{1}_{\{V_n^u>\frac{1}{1+\rho^{-(2B^u_n-n)}}\}}\mathbb{1}_{\{X_n^u=-1\}}\Big|\mathcal{G}_1\right]\\
&=&\sum_{n=n_3}^\infty\mathbb{1}_{\{X_n^u=-1\}}\,\mathbb{E}\left[\mathbb{1}_{\{V_n^u>\frac{1}{1+\rho^{-(2B^u_n-n)}}\}}\Big|\mathcal{G}_1\right].
\end{eqnarray*}

For the two last steps we use the fact that both $n_3$ and $X_n^u$ are $\mathcal{G}_1$ measurable. Then, since $V_n^u$ is independent of $\mathcal{G}_1$ on $\{X_n^u=-1\}$, we have

\begin{eqnarray*}
\mathbb{E}\left[\sum_{n=n_3}^\infty\mathbb{1}_{\mathfrak{c}_n^u={\mathfrak{w}}}\Big|\mathcal{G}_1\right]&\leq&\sum_{n=n_3}^\infty\mathbb{1}_{\{X_n^u=-1\}}\,\mathbb{E}\left[\mathbb{1}_{\{V_n^u>\frac{1}{1+\rho^{-n(2\beta-1)}}\}}\Big|\mathcal{G}_1\right]\\
&\leq&\sum_{n=n_3}^\infty\mathbb{1}_{\{X_n^u=-1\}}\left(1-\frac{1}{1+\rho^{-n(2\beta-1)}}\right)\\
&\leq&\sum_{n=n_3}^\infty\mathbb{1}_{\{X_n^u=-1\}}\frac{1}{1+\rho^{n(2\beta-1)}}\\
&\leq&\sum_{n=n_3}^\infty\frac{1}{1+\rho^{n(2\beta-1)}}\\
&<&\infty.
\end{eqnarray*}

The above calculation implies that the number of white balls in each urn $u$ is finite a.s. This proves the theorem.
\end{proof}

\subsection{The critical case}

In this subsection we consider the case $p=1/2$.
It turns out that in this setting, the {\RUPa} behave as in the supercritical case. Nonetheless, the argument needs some small adaptations. 

In the last argument we use conditioning with respect to $\mathcal{G}_1$ which contains informations about the future. In the following proof we will avoid this kind of conditioning with respect to the future nonetheless, either method could be used in both settings.

We first state another consequence to the law of large numbers that is again left to the reader. (cf. Lemma \ref{lemme1}.)

\begin{lemme}
\label{lemme2}
Fix $v\in\mathcal{U}$. Let $Y_1, Y_2, Y_3, \cdots$ be a sequence of independent random variables with the following law:
\begin{itemize}
	\item if $n\equiv v \mod U$, then $\mathbb{P}[Y_n=1]=1$,
	\item if $n\not\equiv v \mod U$, then $\mathbb{P}[Y_n=1]=1-\mathbb{P}[Y_n=-1]=p =\frac{1}{2}$.
\end{itemize}
Then,
$$\forall \alpha\in(0,1/U), \mathbb{P}\left[\forall n, Y_1+\cdots+Y_n\geq\alpha n\right]>0$$
\end{lemme}

Set $\theta^\prime$ to be the above probability. (Again $\theta^\prime$ depends on the choice of $\alpha$.)

\begin{theorem}
\label{pcritic}
Assume that $p=1/2$.
$$\mathbb{P}\left[\exists n_0\in\N, \exists \mathfrak{d}\in\mathfrak{C} \text{ such that } \forall n\geq n_0, \forall u\in\mathcal{U} \text{ we have } \mathfrak{c}_n^u = \mathfrak{d}\right]=1$$
In words, all the balls drawn after some finite time $n_0$ have the same color $\mathfrak{d}$.
\end{theorem}

For any $n\in\N$ define the stopping time:
$$T^n=\inf\left\{m> n ~:~~ \forall u\in\mathcal{U}, \sum_{i=n+1}^m X_i^u\geq n\right\}.$$

Note that $(\sum_{i=n+1}^m X_i^u, u\in\mathcal{U})_{m>n}$ is a classical random walk on $\Z^U$, and that $T^n$ is the time when this walk enters the domain $[n,\infty)^{U}$ so we know that
$$T^n<\infty ~~\text{a.s.}$$

Before proving Theorem \ref{pcritic}, we state a lemma.

\begin{lemme}\label{lemmecri3}
Assume that $p=1/2$.
There exist a constant $\eta_1>0$ such that for all $n_1\in\N$,
$$\mathbb{P}\left[\bigcap_{n> n_1}\bigcap_{u\in\mathcal{U}}\left\{\{n\leq T^{n_1}\}\cap \{X_n^u=-1\}\right\}\cup \{\mathfrak{c}_n^u=\mathfrak{m}_{n_1}^*\}\Big| \mathcal{F}_{n_1}\right]\geq\eta_1~\text{ a.s.}$$
This means that there is a positive probability uniformly bounded away from 0 so that
\begin{itemize}
\item during $[n_1+1,T^{n_1}]$ each time a ball is drawn out of the $U$ urns combined, it is of the majority color at time $n_1$;
\item after time $T^{n_1}$, all the balls drawn are of this majority color.
\end{itemize}
\end{lemme}

\begin{proof}
Without loss of generality, in this proof we will suppose that the majority color at time $n_1$ is black, in symbols, $\mathfrak{m}_{n_1}^*=\mathfrak{b}$.

We first prove that there is a positive probability uniformly bounded away from 0 that up to $T^{n_1}$, all the balls drawn out of all the $U$ urns combined are black. The second step will be to prove that with a positive probability all the balls drawn after $T^{n_1}$ are black.

\textbf{First step:}
Since $\mathfrak{b}=\mathfrak{m}_{n_1}^*$ is the majority color at $n_1$ there is at least one urn $v\in\mathcal{U}$ such that its majority color is $\mathfrak{m}_{n_1}^*=\mathfrak{b}$, that is $\mathfrak{m}_{n_1}^v=\mathfrak{m}_{n_1}^*=\mathfrak{b}$.

For all $n\in\N$ and $u\in \mathcal{U}$ set $k=nU+u$ and
$$Y_{k}=\mathbb{1}_{\{X^{\prime}_k=1\}\cup\{u=v\}}-\mathbb{1}_{\{X^{\prime}_k=-1\}\cap\{u\neq v\}}.$$
Set $k_1=(n_1+1)U$, the sequence $Y_{k_1+1}, Y_{k_1+2}, Y_{k_1+3},\cdots$ satisfies the conditions of Lemma \ref{lemme2} and is independent of $\mathcal{F}_{n_1}$. So for any $\alpha\in(0,1/U)$ we have:

$$\mathbb{P}\left[A|\mathcal{F}_{n_1}\right]=\mathbb{P}\left[A\right]=\theta^{\prime}>0,$$
where $A= \{Y_{k_1+1}+ \cdots+Y_k\geq\alpha (k-k_1),\forall k\geq k_1\}$. Denote by $C^{\prime}_k=C^u_n$ the event
$$C^{\prime}_k=C^u_n=\left(\{X^{\prime}_k=-1\}\cap\{u\neq v\}\right)\cup \left\{\mathfrak{c}_k^{\prime}=\mathfrak{m}_{n_1}^*\right\}$$
that is \textit{``if the draw for urn $u$ at time $n$ is either out of {\drawall} or out of $v$ then its color is $\mathfrak{m}_{n_1}^*$.''} Then we have:

\begin{eqnarray*}
\mathbb{P}\left[\left.\bigcap_{n=n_1+1}^{T^{n_1}}\bigcap_{u\in\mathcal{U}}\{X_n^u=-1\}\cup \{\mathfrak{c}_n^u=\mathfrak{b}\}\right| \mathcal{F}_{n_1}\right]&=&\mathbb{P}\left[\left.\bigcap_{k=k_1+1}^{(T^{n_1}+1)U}\{X_k^{\prime}=-1\}\cup \left\{\mathfrak{c}_k^{\prime}=\mathfrak{b}\right\}\right| \mathcal{F}_{n_1}\right]\\
&\geq&\mathbb{P}\left[\left.\bigcap_{k=k_1+1}^{(T^{n_1}+1)U}C^{\prime}_k\right| \mathcal{F}_{n_1}\right]\\
&\geq&\mathbb{P}\left[\left.A\cap \bigcap_{k=k_1+1}^{(T^{n_1}+1)U}C^{\prime}_k\right| \mathcal{F}_{n_1}\right]\\
&\geq&\mathbb{E}\left[\left.\mathbb{1}_A\times\prod_{k=k_1+1}^{(T^{n_1}+1)U}\mathbb{1}_{C^{\prime}_k} \right| \mathcal{F}_{n_1}\right]\\
&\geq&\mathbb{E}\left[\left.\mathbb{1}_A\times\mathbb{E}\left[\left.\prod_{k=k_1+1}^{(T^{n_1}+1)U}\mathbb{1}_{C^{\prime}_k}\right| \mathcal{F}_{n_1};A \right]\right| \mathcal{F}_{n_1}\right].
\end{eqnarray*}

We next give a lower bound for the interior expectation. Fix $K> k_1$ and set $\mathcal{F}_{n_1}^K=\sigma(\mathcal{F}_{n_1};A,C^{\prime}_{k_1+1},...,C^{\prime}_{K-1})$. Then, on $A$, we have

\begin{eqnarray*}
\mathbb{E}\left[\left.\prod_{k=k_1+1}^{K}\mathbb{1}_{C^{\prime}_k}\right| \mathcal{F}_{n_1};A \right]&=&\mathbb{E}\left[\left.\mathbb{E}\left[\left.\mathbb{1}_{C^{\prime}_K}\right|\mathcal{F}_{n_1}^K\right]\prod_{k=k_1+1}^{K-1}\mathbb{1}_{C^{\prime}_k}\right| \mathcal{F}_{n_1};A \right],
\end{eqnarray*}
where
\begin{eqnarray*}
\mathbb{E}\left[\mathbb{1}_{C^{\prime}_K}\big|\mathcal{F}_{n_1}^K\right]&=&\mathbb{E}\left[\left.\mathbb{1}_{\left(\{X^{\prime}_K=-1\}\cap\{u\neq v\}\right)\cup \left\{\mathfrak{c}_K^{\prime}=\mathfrak{m}_{n_1}^*\right\}}\right|\mathcal{F}_{n_1}^K\right], \text{ with }u\equiv K \mod U.\\
\end{eqnarray*}

Therefore, on $A$ and $\cap_{k=k_1+1}^{K-1}C^{\prime}_k$, there are three cases to consider.
\begin{enumerate}
\item If $X^{\prime}_K=1$:
\begin{eqnarray*}
\mathbb{E}\left[\mathbb{1}_{C^{\prime}_K}\big|\mathcal{F}_{n_1}^K\right]&=&\mathbb{E}\left[\left.\mathbb{1}_{\left\{\mathfrak{c}_K^{\prime}=\mathfrak{m}_{n_1}^*\right\}}\right|\mathcal{F}_{n_1}^K\right]\\
&=&\mathbb{E}\left[\left.\mathbb{1}_{\left\{V^{\prime}_K<[{1+\rho^{-B_{\lfloor\frac{K-1}{U}\rfloor}^*+W_{\lfloor\frac{K-1}{U}\rfloor}^*}}]^{-1}\right\}}\right|\mathcal{F}_{n_1}^K\right],\\
\end{eqnarray*}
and on $A$ and $\cap_{k=k_1+1}^{K-1}C^{\prime}_k$, we have:
\begin{eqnarray*}
B_{\lfloor\frac{K-1}{U}\rfloor}^*-W_{\lfloor\frac{K-1}{U}\rfloor}^*&=&B_{n_1}^*-W_{n_1}^*+\sum_{n=n_1+1}^{\lfloor\frac{K-1}{U}\rfloor}\sum_{u=1}^U\left(\mathbb{1}_{\{\mathfrak{c}^{u}_n=\mathfrak{b}\}}-\mathbb{1}_{\{\mathfrak{c}^{u}_n=\mathfrak{w}\}}\right)\\
&\geq&\sum_{n=n_1+1}^{\lfloor\frac{K-1}{U}\rfloor}\sum_{u=1}^U\left(\mathbb{1}_{\{\mathfrak{c}^{u}_n=\mathfrak{b}\}}-\mathbb{1}_{\{\mathfrak{c}^{u}_n=\mathfrak{w}\}}\right)\\
&\geq&\sum_{n=n_1+1}^{\lfloor\frac{K-1}{U}\rfloor}\sum_{u=1}^U \left(\mathbb{1}_{\{X^{u}_n=1\}\cup\{u=v\}}-\mathbb{1}_{\{X^{u}_n=-1\}\cap\{u\neq v\}}\right) \text{ on }\cap_{k=k_1+1}^{K-1} C^{\prime}_k\\
&=&\sum_{n=n_1+1}^{\lfloor\frac{K-1}{U}\rfloor}\sum_{u=1}^U Y_{Un+u}\geq\sum_{k=k_1+1}^{k_1+U\left(\lfloor\frac{K-1}{U}\rfloor-n-1)\right)} Y_{k}\\
&\geq&\alpha U(\lfloor\frac{K-1}{U}\rfloor-n_1)\geq\alpha \left(K-1-k_1\right),\\
\end{eqnarray*}

where the second to last inequality is due to $A$. Since $V^{\prime}_K$ is independent of $\mathcal{F}_{n_1}^K$, we deduce that:

\begin{eqnarray*}
\mathbb{E}\left[\mathbb{1}_{C^{\prime}_K}\big|\mathcal{F}_{n_1}^K\right]&\geq&\mathbb{E}\left[\mathbb{1}_{\left\{V^{\prime}_K<\left[1+\rho^{-\alpha\left(K-1-k_1\right)}\right]^{-1}\right\}}\Bigg|\mathcal{F}_{n_1}^K\right]\\
&=&\frac{1}{1+\rho^{-\alpha\left(K-1-k_1\right)}}.\\
\end{eqnarray*}

\item If $X^{\prime}_K=-1$ and $u=v$, then
\begin{eqnarray*}
\mathbb{E}\left[\mathbb{1}_{C^{\prime}_K}\big|\mathcal{F}_{n_1}^K\right]&=&\mathbb{E}\left[\mathbb{1}_{\left\{\mathfrak{c}_K^{\prime}=\mathfrak{m}_{n_1}^*\right\}}\big|\mathcal{F}_{n_1}^K\right]\\
&=&\mathbb{E}\left[\mathbb{1}_{\left\{V^{\prime}_K<\left[{1+\rho^{-\left(B_{K-1}^v-W_{K-1}^v\right)}}\right]^{-1}\right\}}\big|\mathcal{F}_{n_1}^K\right].\\
\end{eqnarray*}

On $A$ and $\cap_{k=k_1+1}^{K-1}C^{\prime}_k$, we have:
\begin{eqnarray*}
B_{K-1}^v-W_{K-1}^v&=&B_{n_1}^v-W_{n_1}^v+\sum_{n=n_1+1}^{\lfloor\frac{K-1}{U}\rfloor}\left(\mathbb{1}_{\{\mathfrak{c}^{v}_n=\mathfrak{b}\}}-\mathbb{1}_{\{\mathfrak{c}^{v}_n=\mathfrak{w}\}}\right)\\
&\geq&\sum_{n=n_1+1}^{\lfloor\frac{K-1}{U}\rfloor}\left(\mathbb{1}_{\{\mathfrak{c}^{v}_n=\mathfrak{b}\}}-\mathbb{1}_{\{\mathfrak{c}^{v}_n=\mathfrak{w}\}}\right)\\
&\geq&\sum_{n=n_1+1}^{\lfloor\frac{K-1}{U}\rfloor}\left(\mathbb{1}_{\{X^v_n=1\}\cup\{u=v\}}-\mathbb{1}_{\{X_n^v=-1\}\cap\{u\neq v\}}\right), \text{ on }\cap_{h=k_1+1}^{K-1} C^{\prime}_k\\
&\geq&\sum_{n=n_1+1}^{\lfloor\frac{K-1}{U}\rfloor}1, \text{ because }u=v\\
&\geq&\lfloor\frac{K-1}{U}\rfloor-n_1 \geq \frac{1}{U}\left(K-1-k_1\right)\\
&\geq&\alpha \left(K-1-k_1\right),\\
\end{eqnarray*}

where in the last line we use the fact $\alpha<\frac{1}{U}$. Since $V^{\prime}_K$ is independent of $\mathcal{F}_{n_1}^K$, we deduce that:

\begin{eqnarray*}
\mathbb{E}\left[\mathbb{1}_{C^{\prime}_K}\big|\mathcal{F}_{n_1}^K\right]&\geq&\mathbb{E}\left[\mathbb{1}_{\left\{V^{\prime}_K<\frac{1}{1+\rho^{-\alpha\left(K-1-k_1\right)}}\right\}}\Bigg|\mathcal{F}_{n_1}^K\right]\\
&\geq&\frac{1}{1+\rho^{-\alpha\left(K-1-k_1\right)}}.\\
\end{eqnarray*}

\item If $X^{\prime}_K=-1$ and $u\neq v$, then
\begin{eqnarray*}
\mathbb{E}\left[\mathbb{1}_{C^{\prime}_K}\big|\mathcal{F}_{n_1}^K\right]&=&\mathbb{E}\left[1\big|\mathcal{F}_{n_1}^K\right]\\
&=&1.
\end{eqnarray*}
\end{enumerate}

In all the three cases we have:
$$\mathbb{E}\left[\mathbb{1}_{C^{\prime}_K}\big|\mathcal{F}_{n_1}^K\right]\geq\frac{1}{1+\rho^{-\alpha\left(K-1-k_1\right)}}.$$

As a conclusion, we have, on the event $A$:

\begin{eqnarray*}
\mathbb{E}\left[\prod_{k=k_1+1}^{K}\mathbb{1}_{C^{\prime}_k}\big| \mathcal{F}_{n_1};A \right]&=&\mathbb{E}\left[\mathbb{E}\left[\mathbb{1}_{C^{\prime}_K}\big|\mathcal{F}_{n_1}^K\right]\prod_{k=k_1+1}^{K-1}\mathbb{1}_{C^{\prime}_k}\big| \mathcal{F}_{n_1};A \right]\\
&\geq&\mathbb{E}\left[\frac{1}{1+\rho^{-\alpha\left(K-1-k_1\right)}}\prod_{k=k_1+1}^{K-1}\mathbb{1}_{C^{\prime}_k}\big| \mathcal{F}_{n_1};A \right]\\
&\geq&\frac{1}{1+\rho^{-\alpha\left(K-1-k_1\right)}}\mathbb{E}\left[\prod_{k=k_1+1}^{K-1}\mathbb{1}_{C^{\prime}_k}\big| \mathcal{F}_{n_1};A \right].
\end{eqnarray*}

And then by induction:

\begin{eqnarray*}
\mathbb{E}\left[\prod_{k=k_1+1}^{K}\mathbb{1}_{C^{\prime}_k}\big| \mathcal{F}_{n_1};A \right]&\geq&\prod_{k=k_1}^{K-1}\frac{1}{1+\rho^{-\alpha\left(k-k_1\right)}}\geq\prod_{k=0}^{K-1-k_1}\frac{1}{1+\rho^{-\alpha k}}\geq\prod_{k=0}^{\infty}\frac{1}{1+\rho^{-\alpha k}}.
\end{eqnarray*}

Now take $K=(T^{n_1}+1)U$, we obtain:
\begin{eqnarray*}
\mathbb{E}\left[\prod_{k=k_1+1}^{(T^{n_1}+1)U}\mathbb{1}_{C^{\prime}_k}\big| \mathcal{F}_{n_1};A \right]&\geq&\prod_{k=0}^{\infty}\frac{1}{1+\rho^{-\alpha k}}.
\end{eqnarray*}
We can now conclude that:
\begin{eqnarray*}
\mathbb{P}\left[\bigcap_{n= n_1+1}^{(T^{n_1}+1)-1}\bigcap_{u\in\mathcal{U}}\{X_n^u=-1\}\cup \{\mathfrak{c}_n^u=\mathfrak{b}\}\Big| \mathcal{F}_{n_1}\right]&\geq&\mathbb{E}\left[\mathbb{1}_A\times\mathbb{E}\left[\prod_{k=k_1+1}^{(T^{n_1}+1)U}\mathbb{1}_{C^{\prime}_k}\Big| \mathcal{F}_{n_1};A \right]\big| \mathcal{F}_{n_1}\right]\\
&\geq&\mathbb{E}\left[\mathbb{1}_A\times\prod_{k=0}^{\infty}\frac{1}{1+\rho^{-\alpha k}}\big| \mathcal{F}_{n_1}\right]\\
&\geq&\prod_{k=0}^{\infty}\frac{1}{1+\rho^{-\alpha k}}\times\mathbb{E}\left[\mathbb{1}_A\big| \mathcal{F}_{n_1}\right]\\
&\geq&\theta^\prime\prod_{k=0}^{\infty}\frac{1}{1+\rho^{-\alpha k}}.
\end{eqnarray*}
This ends the first step of the argument.

\textbf{Second step:}
Denote by $S_1$ the event of the step 1, that is
$$S_1=\bigcap_{n= n_1+1}^{T^{n_1}}\bigcap_{u\in\mathcal{U}}\{X_n^u=-1\}\cup \{\mathfrak{c}_n^u=\mathfrak{b}\}.$$
We now consider the probability that on $S_1$, all the balls drawn after $T^{n_1}$ are black. Let us first prove that on $S_1$ the majority color in each of the $U$ urns at time $T^{n_1}$ is black, that is
$$\forall u\in\mathcal{U},~~\mathfrak{m}_{T^{n_1}}^u=\mathfrak{m}_{n_1}^*=\mathfrak{b}.$$
On $S_1$, for any given $u\in\mathcal{U}$ we have:
\begin{eqnarray*}
B_{T^{n_1}}^u-W_{T^{n_1}}^u&=&\sum_{n=1}^{T^{n_1}}\mathbb{1}_{\{\mathfrak{c}^u_{n}=\mathfrak{b}\}}-\sum_{n=1}^{T^{n_1}}\mathbb{1}_{\{\mathfrak{c}^u_{n}=\mathfrak{w}\}}\\
&=&\sum_{n=1}^{n_1}\mathbb{1}_{\{\mathfrak{c}^u_{n}=\mathfrak{b}\}}+\sum_{n=n_1+1}^{T^{n_1}}\mathbb{1}_{\{\mathfrak{c}^u_{n}=\mathfrak{b}\}}-\sum_{n=1}^{n_1}\mathbb{1}_{\{\mathfrak{c}^u_{n}=\mathfrak{w}\}}-\sum_{n=n_1}^{T^{n_1}}\mathbb{1}_{\{\mathfrak{c}^u_{n}=\mathfrak{w}\}}\\
&\geq&0+\sum_{n=n_1+1}^{T^{n_1}}\mathbb{1}_{\{\mathfrak{c}^u_{n}=\mathfrak{b}\}}-n_1-\sum_{n=n_1}^{T^{n_1}}\mathbb{1}_{\{\mathfrak{c}^u_{n}=\mathfrak{w}\}}\\
&\geq&\sum_{n=n_1+1}^{T^{n_1}}\mathbb{1}_{\{X_n^u=1\}}-\sum_{n=n_1}^{T^{n_1}}\mathbb{1}_{\{X_n^{u}=-1\}}-n_1,\text{   because we are on }S_1\\
&\geq&\sum_{n=n_1+1}^{T^{n_1}}X_n^u-n_1\\
&\geq&n_1-n_1\text{  by the definition of }T^{n_1}\\
&\geq&0.
\end{eqnarray*}
This proves that $\mathfrak{m}_{T^{n_1}}^u=\mathfrak{b}$. Let us now compute the probability that all the balls drawn after $T^{n_1}$ are black. Take $N>T^{n_1}+1$ then
\begin{eqnarray*}
&&\mathbb{P}\left[\left.\bigcap_{n=T^{n_1}+1}^{T^{n_1}+N}\bigcap_{u\in\mathcal{U}}\{\mathfrak{c}_n^u=\mathfrak{b}\}\right|\mathcal{F}_{n_1};S_1\right]\\
&&~~~~~~~~~~~~~~=\mathbb{E}\left[\left.\prod_{n=T^{n_1}+1}^{T^{n_1}+N}\prod_{u\in\mathcal{U}}\mathbb{1}_{\{\mathfrak{c}_n^u=\mathfrak{b}\}}\right|\mathcal{F}_{n_1};S_1\right]\\
&&~~~~~~~~~~~~~~=\mathbb{E}\left[\left.\prod_{n=T^{n_1}+1}^{T^{n_1}+N-1}\prod_{u\in\mathcal{U}}\mathbb{1}_{\{\mathfrak{c}_n^u=\mathfrak{b}\}}\times\mathbb{E}\left[\left.\prod_{u\in\mathcal{U}}\mathbb{1}_{\{\mathfrak{c}_{T^{n_1}+N}^u=\mathfrak{b}\}}\right|\mathcal{F}_{T^{n_1}+N-1};S_1\right]\right|\mathcal{F}_{n_1};S_1\right]\\
&&~~~~~~~~~~~~~~\geq\mathbb{E}\left[\left.\prod_{n=T^{n_1}+1}^{T^{n_1}+N-1}\prod_{u\in\mathcal{U}}\mathbb{1}_{\{\mathfrak{c}_n^u=\mathfrak{b}\}}\times\prod_{u\in\mathcal{U}}\frac{1}{1+\rho^{W_{T^{n_1}+N-1}^u-B_{T^{n_1}+N-1}^u}}\right|\mathcal{F}_{n_1};S_1\right].\\
\end{eqnarray*}
On $S_1 \cap \bigcap_{n=T^{n_1}+1}^{T^{n_1}+N-1}\bigcap_{u\in\mathcal{U}}\{\mathfrak{c}_n^u=\mathfrak{b}\}$, we have
$$B_{T^{n_1}+N-1}^u-W_{T^{n_1}+N-1}^u>N-1.$$
Therefore,
\begin{eqnarray*}
\mathbb{P}\left[\left.\bigcap_{n=T^{n_1}+1}^{T^{n_1}+N}\bigcap_{u\in\mathcal{U}}\{\mathfrak{c}_n^u=\mathfrak{b}\}\right|\mathcal{F}_{n_1};S_1\right]&\geq&\mathbb{E}\left[\left.\prod_{n=T^{n_1}+1}^{T^{n_1}+N-1}\prod_{u\in\mathcal{U}}\mathbb{1}_{\{\mathfrak{c}_n^u=\mathfrak{b}\}}\times\prod_{u\in\mathcal{U}}\frac{1}{1+\rho^{-(N-1)}}\right|\mathcal{F}_{n_1};S_1\right]\\
&\geq&\left(\frac{1}{1+\rho^{-(N-1)}}\right)^U\times\mathbb{E}\left[\left.\prod_{n=T^{n_1}+1}^{T^{n_1}+N-1}\prod_{u\in\mathcal{U}}\mathbb{1}_{\{\mathfrak{c}_n^u=\mathfrak{b}\}}\right|\mathcal{F}_{n_1};S_1\right].\\
\end{eqnarray*}
And then by induction :
\begin{eqnarray*}
\mathbb{P}\left[\left.\bigcap_{n=T^{n_1}+1}^{T^{n_1}+N}\bigcap_{u\in\mathcal{U}}\{\mathfrak{c}_n^u=\mathfrak{b}\}\right|\mathcal{F}_{n_1};S_1\right]&\geq&\left(\prod_{n=0}^{N-1}\frac{1}{1+\rho^{-n}}\right)^U\\
&\geq&\left(\prod_{n=0}^{\infty}\frac{1}{1+\rho^{-n}}\right)^U.\\
\end{eqnarray*}
Letting $N$ going to the infinity gives us the conclusion of the second step. Set
$$F_{n_1}=\bigcap_{n> n_1}\bigcap_{u\in\mathcal{U}}\left\{\{n\leq T^{n_1}\}\cap \{X_n^u=-1\}\right\}\cup \{\mathfrak{c}_n^u=\mathfrak{b}\}.$$
Therefore we have:
\begin{eqnarray*}
\mathbb{P}\left[\left.A\right| \mathcal{F}_{n_1}\right]
&=&\mathbb{E}\left[\left.\prod_{n= n_1+1}^{\infty}\prod_{u\in\mathcal{U}}\mathbb{1}_{\left\{\{n\leq T^{n_1}\}\cap\{X_n^u=-1\}\right\}\cup \{\mathfrak{c}_n^u=\mathfrak{b}\}}\right| \mathcal{F}_{n_1}\right]\\
&=&\mathbb{E}\left[\left.\prod_{n= n_1+1}^{T^{n_1}}\prod_{u\in\mathcal{U}}\mathbb{1}_{\{X_n^u=-1\}\cup \{\mathfrak{c}_n^u=\mathfrak{b}\}}\times\mathbb{E}\left[\left.\prod_{n= T^{n_1}+1}^{\infty}\prod_{u\in\mathcal{U}}\mathbb{1}_{\{\mathfrak{c}_n^u=\mathfrak{b}\}}\right| \mathcal{F}_{n_1};S_1\right]\right| \mathcal{F}_{n_1}\right]\\
&\geq&\mathbb{E}\left[\left.\prod_{n= n_1+1}^{T^{n_1}}\prod_{u\in\mathcal{U}}\mathbb{1}_{\{X_n^u=-1\}\cup \{\mathfrak{c}_n^u=\mathfrak{b}\}}\times\left(\prod_{n=0}^{\infty}\frac{1}{1+\rho^{-n}}\right)^U\right| \mathcal{F}_{n_1}\right]\text{ by step 2,}\\
&\geq&\theta^\prime\prod_{k=0}^{\infty}\frac{1}{1+\rho^{-\alpha k}}\times\left(\prod_{n=0}^{\infty}\frac{1}{1+\rho^{-n}}\right)^U\text{ by step 1.}\\
\end{eqnarray*}
Hence we obtain Lemma \ref{lemmecri3} with
$$\eta_1=\theta^\prime\prod_{k=0}^{\infty}\frac{1}{1+\rho^{-\alpha k}}\times\left(\prod_{n=0}^{\infty}\frac{1}{1+\rho^{-n}}\right)^U>0.$$
\end{proof}

\begin{proof} (of Theorem \ref{pcritic})
Set
$$F=\bigcup_{\mathfrak{c}\in\{\mathfrak{b},\mathfrak{w}\}} \bigcup_{n_1\in\N} \bigcap_{n> n_1}\bigcap_{u\in\mathcal{U}}\left(\{n\leq T^{n_1}\}\cap \{X_n^u=-1\}\right)\cup \{\mathfrak{c}_n^u=\mathfrak{c}\}.$$

Then:
$$\mathbb{1}_F=\lim_{n_1\rightarrow\infty}\mathbb{P}\left[F|\mathcal{F}_{n_1}\right]\geq\lim_{n_1\rightarrow\infty}\mathbb{P}\left[\left.F_{n_1}\right|\mathcal{F}_{n_1}\right]\geq\eta_1>0~\text{a.s.}$$
by Lemma \ref{lemmecri3}. Therefore $\mathbb{P}(F)=1$. So almost surely, there exits $n_1$ such that
$$\forall n>T^{n_1}, \forall u\in\mathcal{U}, ~~\mathfrak{c}_n^u=\mathfrak{m}_{n_1}^*.$$
To obtain Theorem \ref{pcritic}, just set $\mathfrak{d}=\mathfrak{m}_{n_1}^*$ and $n_0=T^{n_1}+1$.
\end{proof}

\subsection{The subcritical case}

In this subsection we consider the case when $p<1/2$. Note that if $p=0$ then the urns are independent, so each of them chooses a color with probability $1/2$ and draws only this color after a finite time. Hence in this part we will only study the interesting case $0<p<1/2$. One can still remark that Theorem \ref{thsub} is true for $p=0$ as a degenerated case.

\begin{defi}
\label{mero}
Given $x\in\chi$ an environment, an {\up} is said to be $x$-suitable if there exist a color $\mathfrak{d}$ and a finite time $n_0$ such that for each $u\in\mathcal{U}$:
\begin{itemize}
\item either for all $n\geq n_0$, $\mathfrak{c}_u^n=\mathfrak{d}$,
\item or for all $n\geq n_0$, $\mathfrak{c}_u^n=\mathfrak{d}$ each time $x_n^u=1$ and $\mathfrak{c}_u^n\neq \mathfrak{d}$ each time $x_n^u=-1$.
\end{itemize}
We call {\em conformist urns} the urns which satisfy the first condition and {\em nonconformist urns} the urns which satisfy the second condition.
Then denote by $N_c$ the number of conformist urns and by $N_{\not c}=U-N_c$ the number of nonconformist urns.
\end{defi}

If $x$ draw as a sequence of independent Bernoulli variables of parameter $p$, a $x$-suitable {\up} can be schematically represented as follows:
	\begin{center}
	\pspicture(0,-1)(8.5,2.5)

			\psline[border=0pt](0,2)(0,0)
			\psline[border=0pt](0,0)(1,0)
			\psline[border=0pt](1,0)(1,2)
			\psline[border=0pt](0,.6)(1,.6)
			\psline[border=0pt](0,1.8)(1,1.8)
			
			\psline[border=0pt]{-}(1.5,2)(1.5,0)
			\psline[border=0pt]{-}(1.5,0)(2.5,0)
			\psline[border=0pt]{-}(2.5,0)(2.5,2)
			\psline[border=0pt](1.5,.6)(2.5,.6)
			\psline[border=0pt](1.5,1.8)(2.5,1.8)

			\psframe[linecolor=black,fillcolor=black,fillstyle=solid](3,.6)(4,1.8)
			\psline[border=0pt]{-}(3,2)(3,0)
			\psline[border=0pt]{-}(3,0)(4,0)
			\psline[border=0pt]{-}(4,0)(4,2)
			
			\psframe[linecolor=black,fillcolor=black,fillstyle=solid](4.5,.6)(5.5,1.8)
			\psline[border=0pt]{-}(4.5,2)(4.5,0)
			\psline[border=0pt]{-}(4.5,0)(5.5,0)
			\psline[border=0pt]{-}(5.5,0)(5.5,2)

			\psframe[linecolor=black,fillcolor=black,fillstyle=solid](6,.6)(7,1.8)
			\psline[border=0pt]{-}(6,2)(6,0)
			\psline[border=0pt]{-}(6,0)(7,0)
			\psline[border=0pt]{-}(7,0)(7,2)
			
			\psframe[linecolor=black,fillcolor=black,fillstyle=solid](7.5,.6)(8.5,1.8)
			\psline[border=0pt]{-}(7.5,2)(7.5,0)
			\psline[border=0pt]{-}(7.5,0)(8.5,0)
			\psline[border=0pt]{-}(8.5,0)(8.5,2)
			
			\rput(1.25,-.2){$\underbrace{~~~~~~~~~~~~~~~~~~}$}
			\rput(5.75,-.2){$\underbrace{~~~~~~~~~~~~~~~~~~~~~~~~~~~~~~~~~~~~~~~~~}$}
			\rput(1.25,-.5){{conformist urns}}
			\rput(5.75,-.5){{nonconformist urns}}
			\rput(-.1,1.3){$\left\{\raisebox{0pt}[20pt][0pt]{~}\right.$}
			\rput(-.1,.3){$\left\{\raisebox{0pt}[7.5pt][0pt]{~}\right.$}
			\rput(-.8,1.2){$1-p$}
			\rput(-.5,.3){$p$}

	\endpspicture
	
\noindent
{\footnotesize Figure 6:
Schematic representation of an $x$-suitable {\up}.}
\setcounter{figure}{6}

\vspace{0.3cm}
\noindent
	\end{center}

\begin{defi}(Conformism Equation)
Suppose that for some $x\in\chi$, the {$\RUPa^x$} is $x$-suitable almost surely, then it is said to satisfy the Conformism Equation if:
\begin{equation}
\label{CE}
\tag{CE}
(1-p)N_{\not c}<\frac{U}{2}.
\end{equation}
\end{defi}

Note that if $x$ is drawn as a sequence of independent Bernoulli variables of parameter $p$, then almost surely, the asymptotic proportion of the color $\mathfrak{d}$ from Definition \ref{mero} is
$$\lim_{n\to\infty}\frac{D_n^*}{Un}=\frac{N_cn + pN_{\not c}n}{Un}=\frac{(U-N_{\not c})+pN_{\not c}}{U}=1-\frac{(1-p)N_{\not c}}{U}>\frac{1}{2}\text{  under (\ref{CE}}),$$
where $D_n^*$ is the total number of balls of color $\mathfrak{d}$ in the $U$ urns combined at time $n$.

\begin{prop}
\label{propsub}
For almost every environment $x$, if $r$ is an $x$-suitable {\up} which satisfies (\ref{CE}), then
$$P^x[r]>0.$$
\end{prop}

\begin{proof}
Let $x$ be an environment that is sampled as a sequence $(x_n^u)_{n\in\N,u\in\mathcal{U}}$ of i.i.d random variables such that
$$\mathbb{P}[x_n^u=1]=1-\mathbb{P}[x_n^u=-1]=p, \forall n\in\N,\forall u\in\mathcal{U}.$$
Then $x$ almost surely satisfies:
\begin{equation}
\forall u\in\mathcal{U}, \lim_{n\to\infty}\frac{\sum_{i=1}^{n}\mathbb{1}_{\{x_n^u=1\}}}{n}=p,
\end{equation}
or equivalently
\begin{equation}
\forall u\in\mathcal{U}, \lim_{n\to\infty}\frac{\sum_{i=1}^{n}\mathbb{1}_{\{x_n^u=-1\}}}{n}=1-p.
\end{equation}
Now, let $r$ be a deterministic $x$-suitable {\RUPa} which satisfies (\ref{CE}) and set $\mathfrak{d}$ and $n_0$ as in Definition \ref{mero}, and denote by $d_n^u$ (resp. $d_n^*$) the number of balls of color $\mathfrak{d}$ in urn $u$ at time $n$ (resp. in all the urns at time $n$). Then
$$\lim_{n\rightarrow\infty}\frac{d_n^*}{n}=\frac{N_c+p N_{\not c}}{U}=\frac{U-N_{\not c}+p N_{\not c}}{U}=1-\frac{(1-p)N_{\not c}}{U}>\frac{1}{2},$$
due to (\ref{CE}). For each conformist urn $u$,
$$\lim_{n\rightarrow\infty}\frac{d_n^u}{n}=1,$$
and for each nonconformist urn $u$:
$$\lim_{n\rightarrow\infty}\frac{n-d_n^u}{n}=1-p>\frac{1}{2}.$$

Fix $\beta\in(1/2,\min \{1-p,1-(1-p)N_{\not c}/U\})$. There exist a time $n_1>n_0$ such that for all $n\geq n_1$:
\begin{itemize}
\item $\frac{d_n^*}{Un}>\beta$,
\item $\frac{d_n^u}{n}>\beta$ for each conformist urn $u$,
\item $\frac{n-d_n^u}{n}>\beta$ for each nonconformist urn $u$.
\end{itemize}
Let $R$ be an $\RUPa^x$. Then
\begin{eqnarray*}
P^x\left[R=r\right]&=&E^x\left[\mathbb{1}_{\{R=r\}}\right],\\
&=&E^x\left[E^x\left[\left.\mathbb{1}_{\{(R_n)_{n>n_1}=(r_n)_{n>n_1}\}}\right|\mathcal{F}_{n_1}\right]\mathbb{1}_{\{(R_n)_{0\leq n\leq n_1}=(r_n)_{0\leq n\leq n_1}\}}\right].\\
\end{eqnarray*}
Fix $N> n_1$ then on $\{(R_n)_{0\leq n\leq n_1}=(r_n)_{0\leq n\leq n_1}\}$ we have:
\begin{eqnarray*}
E^x\left[\left.\prod_{n=n_1+1}^{N}\mathbb{1}_{\{R_n=r_n\}}\right|\mathcal{F}_{n_1}\right]&=&E^x\left[\left.E^x\left[\left.\mathbb{1}_{\{R_N=r_N\}}\right|\mathcal{F}_{N-1}\right]\prod_{n=n_1+1}^{N-1}\mathbb{1}_{\{R_n=r_n\}}\right|\mathcal{F}_{n_1}\right].
\end{eqnarray*}
And on $\{(R_n)_{0\leq n\leq N-1}=(r_n)_{0\leq n\leq N-1}\}$ we have
\begin{eqnarray*}
E^x\left[\left.\mathbb{1}_{\{R_N=r_N\}}\right|\mathcal{F}_{N-1}\right]&=&E^x\left[\left.\prod_{u\in\mathcal{U}}\mathbb{1}_{\{R_N(u)=r_N(u)\}}\right|\mathcal{F}_{N-1}\right]\\
&=&\prod_{u\in\mathcal{U}}E^x\left[\left.\mathbb{1}_{\{R_N(u)=r_N(u)\}}\right|\mathcal{F}_{N-1}\right],
\end{eqnarray*}
where
\[
  E^x\left[\left.\mathbb{1}_{\{R_N(u)=r_N(u)\}}\right|\mathcal{F}_{N-1}\right]= \left\{
          \begin{array}{ll}
            \frac{1}{1+\rho^{U{(N-1)}-2d_{N-1}^*}} & \qquad \text{if }x_N^u=1, \\
            \frac{1}{1+\rho^{{N-1}-2d_{N-1}^u}} & \qquad \text{if }x_N^u=-1\text{ and }u\text{ is conformist,} \\
            \frac{1}{1+\rho^{-{(N-1)}+2d_{N-1}^u}} & \qquad \text{if }x_N^u=-1\text{ and }u\text{ is nonconformist.}\\
          \end{array}
        \right.
\]
Therefore
\begin{eqnarray*}
E^x\left[\left.\mathbb{1}_{\{R_N(u)=r_N(u)\}}\right|\mathcal{F}_{N-1}\right]&\geq&\min\left(\frac{1}{1+\rho^{U{(N-1)}-2d_{N-1}^*}},\frac{1}{1+\rho^{{N-1}-2d_{N-1}^u}},\frac{1}{1+\rho^{-{(N-1)}+2d_{N-1}^u}}\right)\\
&\geq&\min\left(\frac{1}{1+\rho^{U{(N-1)}(1-2\beta)}},\frac{1}{1+\rho^{{(N-1)}(1-2\beta)}},\frac{1}{1+\rho^{{(N-1)}(1-2\beta)}}\right)\\
&\geq&\frac{1}{1+\rho^{(N-1)(1-2\beta)}}.
\end{eqnarray*}

By induction it follows that:

$$E^x\left[\left.\prod_{n=n_1+1}^{N}\mathbb{1}_{\{R_n=r_n\}}\right|\mathcal{F}_{n_1}\right]\geq\prod_{n=n_1+1}^N\left(\frac{1}{1+\rho^{{(n-1)}(1-2\beta)}}\right)^U\geq\prod_{n=0}^{\infty}\left(\frac{1}{1+\rho^{n(1-2\beta)}}\right)^U.$$

Let $N$ go to infinity.
\begin{eqnarray*}
E^x\left[\left.\prod_{n=n_1+1}^{\infty}\mathbb{1}_{\{R_n=r_n\}}\right|\mathcal{F}_{n_1}\right]&\geq&\prod_{n=0}^{\infty}\left(\frac{1}{1+\rho^{n(1-2\beta)}}\right)^U>0.
\end{eqnarray*}

Set $\eta=P^x\left[(R_n)_{0\leq n\leq n_1}=(r_n)_{0\leq n\leq n_1}\right]>0$ then we have:
\begin{eqnarray*}
P^x\left[R=r\right]&\geq&\eta\prod_{n=0}^{\infty}\left(\frac{1}{1+\rho^{n(1-2\beta)}}\right)^U>0,
\end{eqnarray*}
which is positive and this proves the proposition.
\end{proof}

\begin{theorem}
\label{thsub}
Assume that $p<1/2$, then for $\mathbb{P}-a.e.~~ x$ sampled as a sequence $(x_n^u)_{n\in\N,u\in\mathcal{U}}$ of i.i.d random variables such that $\mathbb{P}[x_n^u=1]=1-\mathbb{P}[x_n^u=-1]=p$ we have:
$$P^x[\text{The {$\RUPa^x$} is }x\text{-suitable and satisfies (\ref{CE})}]=1.$$
\end{theorem}

\noindent\textbf{Remark:} Note that in the critical and supercritical cases we also obtain $x$-suitables \RUPa s.
They satisfy \ref{CE} since there are no nonconformist urns. But in those cases, Proposition \ref{propsub} is false since $x$-suitable {\ups} with $N_{\not c}>0$ do not occur a.s.

Before we prove the theorem, let us state a lemma:

\begin{lemme}
\label{lemmesous1}
Assume that $p< 1/2$ and fix $u\in\mathcal{U}$. There exist a constant $\eta_1>0$ such that for all $n_1\in\N$,
$$\mathbb{P}\left[\mathbb{P}\left[\bigcap_{n>n_1}\{X_n^u=1\}\cup \{\mathfrak{c}_n^u=\mathfrak{m}_{n_1}^u\}\Big|\mathcal{G}_{-1}\right]=1 \Big| \mathcal{F}_{n_1}\right]\geq\eta_1~~a.s.$$
This means that with a probability uniformly bounded away from $0$, the information $\mathcal{G}_{-1}$ is sufficient to be sure that after time $n_1$, if a ball is drawn out of the urn $u$ alone, its color is the majority color of $u$ at time $n_1$.
\end{lemme}

This lemma is to be compared with Lemma \ref{lemmesur1}. The argument is almost identical, just replacing $p$ by $1-p$, $\mathcal{G}_1$ by $\mathcal{G}_{-1}$, and the balls drawn out of {\drawall} by the balls drawn out of urn $u$ alone. The proof is left to the reader.

We can now deal with the proof of Theorem \ref{thsub}.

\begin{proof}
For each $u$ in $\mathcal{U}$, denote by $F_{n_1}^u$ the event:
$$F_{n_1}^u=\left\{\mathbb{P}\left[\bigcap_{n>n_1}\{X_n^u=1\}\cup \{\mathfrak{c}_n^u=\mathfrak{m}_{n_1}^u\}\Big|\mathcal{G}_{-1}\right]=1\right\}.$$

Then $(F_{n_1}^u)_{n_1\in\N}$ is an increasing sequence of events, and let $F^u$ be its limit
$$F^u=\lim_{n_1\to\infty}F_{n_1}^u=\bigcup_{n_1}F_{n_1}^u.$$
Thus
$$\mathbb{1}_{F^u}=\lim_{n_1\rightarrow\infty}\mathbb{E}\left[\mathbb{1}_{F^u}|\mathcal{F}_{n_1}\right]=\lim_{n_1\rightarrow\infty}\mathbb{P}\left[F^u|\mathcal{F}_{n_1}\right]\geq\lim_{n_1\rightarrow\infty}\mathbb{P}\left[F_{n_1}^u|\mathcal{F}_{n_1}\right]\geq\eta_1>0\text{ a.s.},$$
by Lemma \ref{lemmesous1}. Therefore $\mathbb{P}\left[F^u\right]=1$ and $\mathbb{P}\left[\bigcap_{u\in\mathcal{U}}F^u\right]=1$.

Denote by $n_2$ the finite random number defined by:
$$n_2=\min\left\{n_1\in\N:\bigcap_{u\in\mathcal{U}}F_{n_1}^u\text{ happens}\right\}.$$

Note that $n_2$ is $\mathcal{G}_{-1}$-measurable so it depends on a part of the future but is still independent of $V^u_n$ when $X^u_n=1$. This is the important point for the following argument to work. To summarize the argument up to this point, we now know that at time $n_2$, each urn $u$ has chosen a color $\mathfrak{d}^u$ such that each time a ball is drawn out of $u$ alone after time $n_2$ then it is almost surely of this color.

Denote by $B$ the number of urns that choose black ($\mathfrak{b}$): $B= \sum_{u\in\mathcal{U}}\mathbb{1}_{\{\mathfrak{d}^u=\mathfrak{b}\}}$ ; and by $W$ the number of urns that choose white ($\mathfrak{w}$): $W= \sum_{u\in\mathcal{U}}\mathbb{1}_{\{\mathfrak{d}^u=\mathfrak{w}\}}$.

Suppose without loss of generality that there is more urns that choose black than urns that choose white, in symbols: $B\geq W$. The fact that in case of ties we choose black has no importance, as will be seen at the end of the argument.

Now we need to distinguish two cases for which the end of the proof will be different. Let us first give in a few lines an idea of what happens in each of the cases:

\textbf{First case:} if $(1-p)B>U/2$ it means that the majority color in {\drawall} will eventually be black. Urns which choose black shall be conformist urns and those which choose white will be nonconformist:

$$N_c=B\text{  and }N_{\not c} = W.$$

\textbf{Second case:} if $(1-p)B<U/2$ (and $(1-p)W<U/2$ since $B\geq W$) then the conformist color is not yet decided, and could be either black or white.
We will see at the end that the equality case $(1-p)B=U/2$ easily follows from the proof of the second case.
\\

Let us start with the proof of the first case. Fix $\beta\in (1/2,(1-p)B/U)$ and denote by $n_3$ the $\mathcal{G}_{-1}$-measurable time:

$$n_3=\inf\left\{m>n_2 : \forall n\geq m, \frac{\sum_{l=n_2}^{n}\sum_{u\in\mathcal{U}}\mathbb{1}_{\{\mathfrak{d}^u=\mathfrak{b}\}\cap\{x_l^u=-1\}}}{nU}\geq\beta\right\}.$$

This time is finite a.s. since

$$\lim_{n\to\infty}\frac{\sum_{l=n_2}^{n}\sum_{u\in\mathcal{U}}\mathbb{1}_{\{\mathfrak{d}^u=\mathfrak{b}\}\cap\{x_l^u=-1\}}}{nU}=(1-p)B/U>\beta.$$

Then for any $n\geq n_3$ we have:

$$\frac{B_n^*}{Un}=\frac{\sum_{l=1}^{n}\sum_{u\in\mathcal{U}}\mathbb{1}_{\{\mathfrak{c}_n^u=\mathfrak{b}\}}}{Un}\geq\frac{\sum_{l=n_2}^{n}\sum_{u\in\mathcal{U}}\mathbb{1}_{\{\mathfrak{d}^u=\mathfrak{b}\}\cap\{x_l^u=-1\}}}{Un}\geq \beta,$$

and therefore for any $N\geq 1$:
\begin{equation*}
\begin{split}
P^x\left[\left.\bigcap_{l=n}^{n+N}\bigcap_{u\in\mathcal{U}}\left\{x_l^u=-1\right\}\cup\left\{\mathfrak{c}_l^u=\mathfrak{b}\right\}\right|\mathcal{G}_{-1},\mathcal{F}_{n}\right]=E^x\left[\left.\prod_{l=n}^{n+N}\prod_{u\in\mathcal{U}}\mathbb{1}_{\left\{x_l^u=-1\right\}\cup\left\{\mathfrak{c}_l^u=\mathfrak{b}\right\}}\right|\mathcal{G}_{-1},\mathcal{F}_{n}\right]\\
=E^x\left[\left.E^x\left[\left.\prod_{u\in\mathcal{U}}\mathbb{1}_{\left\{x_l^u=-1\right\}\cup\left\{\mathfrak{c}_l^u=\mathfrak{b}\right\}}\right|\mathcal{G}_{-1},\mathcal{F}_{n+N-1}\right]\times\prod_{l=n}^{n+N-1}\prod_{u\in\mathcal{U}}\mathbb{1}_{\left\{x_l^u=-1\right\}\cup\left\{\mathfrak{c}_l^u=\mathfrak{b}\right\}}\right|\mathcal{G}_{-1},\mathcal{F}_{n}\right].
\end{split}
\end{equation*}

And on $\bigcap_{l=n}^{n+N-1}\bigcap_{u\in\mathcal{U}}\left\{x_l^u=-1\right\}\cup\left\{\mathfrak{c}_l^u=\mathfrak{b}\right\}$:

\begin{eqnarray*}
E^x\left[\left.\prod_{u\in\mathcal{U}}\mathbb{1}_{\left\{x_l^u=-1\right\}\cup\left\{\mathfrak{c}_l^u=\mathfrak{b}\right\}}\right|\mathcal{G}_{-1},\mathcal{F}_{n+N-1}\right]&=&\prod_{u\in\mathcal{U}}E^x\left[\left.\mathbb{1}_{\left\{x_l^u=-1\right\}\cup\left\{\mathfrak{c}_l^u=\mathfrak{b}\right\}}\right|\mathcal{G}_{-1},\mathcal{F}_{n+N-1}\right]\\
&=&\prod_{u\in\mathcal{U}}E^x\left[\left.\mathbb{1}_{\left\{x_l^u=-1\right\}}+\mathbb{1}_{\left\{x_l^u=1\right\}\cup\left\{\mathfrak{c}_l^u=\mathfrak{b}\right\}}\right|\mathcal{G}_{-1},\mathcal{F}_{n+N-1}\right]\\
&=&\prod_{u\in\mathcal{U}}\mathbb{1}_{\left\{x_l^u=-1\right\}}+\mathbb{1}_{\left\{x_l^u=1\right\}}E^x\left[\left.\mathbb{1}_{\left\{\mathfrak{c}_l^u=\mathfrak{b}\right\}}\right|\mathcal{G}_{-1},\mathcal{F}_{n+N-1}\right]\\
&=&\prod_{u\in\mathcal{U}}\mathbb{1}_{\left\{x_l^u=-1\right\}}+\mathbb{1}_{\left\{x_l^u=1\right\}}\frac{1}{1+\rho^{Ul-2B_l^*}}\\
&=&\prod_{u\in\mathcal{U}}\mathbb{1}_{\left\{x_l^u=-1\right\}}+\mathbb{1}_{\left\{x_l^u=1\right\}}\frac{1}{1+\rho^{Ul(1-2B_l^*/Ul)}}\\
&\geq&\prod_{u\in\mathcal{U}}\mathbb{1}_{\left\{x_l^u=-1\right\}}+\mathbb{1}_{\left\{x_l^u=1\right\}}\frac{1}{1+\rho^{Un(1-2\beta)}}.\\
\end{eqnarray*}

Since $l\geq n \geq n_3$, we have $B_l^*/Ul\geq\beta$. Therefore

\begin{eqnarray*}
E^x\left[\left.\prod_{u\in\mathcal{U}}\mathbb{1}_{\left\{x_l^u=-1\right\}\cup\left\{\mathfrak{c}_l^u=\mathfrak{b}\right\}}\right|\mathcal{G}_{-1},\mathcal{F}_{n+N-1}\right]
&\geq&\prod_{u\in\mathcal{U}}\left(\mathbb{1}_{\left\{x_l^u=-1\right\}}+\mathbb{1}_{\left\{x_l^u=1\right\}}\right)\frac{1}{1+\rho^{Ul(1-2\beta)}}\\
&=&\prod_{u\in\mathcal{U}}\frac{1}{1+\rho^{Ul(1-2\beta)}}\\
&=&\left(\frac{1}{1+\rho^{Ul(1-2\beta)}}\right)^U.\\
\end{eqnarray*}

By induction it follows that

$$P^x\left[\left.\bigcap_{l=n}^{n+N}\bigcap_{u\in\mathcal{U}}\left\{x_l^u=-1\right\}\cup\left\{\mathfrak{c}_l^u=\mathfrak{b}\right\}\right|\mathcal{G}_{-1},\mathcal{F}_{n}\right]\geq \prod_{l=n}^{n+N}\left(\frac{1}{1+\rho^{Ul(1-2\beta)}}\right)^U\geq \prod_{l=0}^{\infty}\left(\frac{1}{1+\rho^{Ul(1-2\beta)}}\right)^U>0.$$

And letting $N\to\infty$, we obtain:

$$P^x\left[\left.\bigcap_{l=n}^{\infty}\bigcap_{u\in\mathcal{U}}\left\{x_l^u=-1\right\}\cup\left\{\mathfrak{c}_l^u=\mathfrak{b}\right\}\right|\mathcal{G}_{-1},\mathcal{F}_{n}\right]\geq\prod_{l=0}^{\infty}\left(\frac{1}{1+\rho^{Ul(1-2\beta)}}\right)^U:=\eta_2>0.$$

Set 
$$F_n=\bigcap_{l=n}^{\infty}\bigcap_{u\in\mathcal{U}}\left\{x_l^u=-1\right\}\cup\left\{\mathfrak{c}_l^u=\mathfrak{b}\right\},$$
and
$$F=\bigcup_{n=n_3}^{\infty}F_n=\bigcup_{n=n_3}^{\infty}\bigcap_{l=n}^{\infty}\bigcap_{u\in\mathcal{U}}\left\{x_l^u=-1\right\}\cup\left\{\mathfrak{c}_l^u=\mathfrak{b}\right\},$$

then by the classical argument we conclude:
$$\mathbb{1}_F=\lim_{n\to\infty}\mathbb{E}\left[\left.\mathbb{E}\left[\left.\mathbb{1}_F\right|\mathcal{F}_n;\mathcal{G}_{-1}\right]\right|\mathcal{F}_n\right]\geq \lim_{n\to\infty}\mathbb{E}\left[\left.\mathbb{E}\left[\left.\mathbb{1}_{F_n}\right|\mathcal{F}_n;\mathcal{G}_{-1}\right]\right|\mathcal{F}_n\right]\geq \eta_2 >0,$$

so that $F$ is an event of full probability, which ends the proof of the first case.
\\

Let us now argue in the second case: suppose that $(1-p)W\leq(1-p)B<U/2$.

For this argument, we need the following lemma:

\begin{lemme}
\label{lemme4}
Fix $\mathcal{E}\subset\mathcal{U}$ a set of urn labels whose cardinality $E:=\left|\mathcal{E}\right|$ satisfies $(1-p)E<U/2$. Let $Y_1, Y_2, Y_3, \cdots$ be a sequence of independent random variables with the following law:
\begin{itemize}
	\item if $n \mod U\in \mathcal{E}$, then $\mathbb{P}[Y_n=1]=1-\mathbb{P}[Y_n=-1]=p <1/2$,
	\item if $n \mod U\not\in \mathcal{E}$, then $\mathbb{P}[Y_n=1]=1$.
\end{itemize}
Then,
$$\forall \alpha\in(0,1-2E(1-p)/U), \mathbb{P}\left[\forall n, Y_1+\cdots+Y_n\geq\alpha n\right]>0.$$
Set $\theta_p^{\prime\prime}$ to be the above probability. (Note that $\theta_p^{\prime\prime}$ depends on the choice of $\alpha$.)
\end{lemme}

Note that this lemma is also a consequence of the law of large numbers which is to be compared with Lemma \ref{lemme1} (supercritical case) and Lemma \ref{lemme2} (critical case). The proof is left to the reader.

Fix $n\geq n_2$ and denote by $\mathfrak{d}=\mathfrak{m}_n^*$ the majority color in {\drawall} at time $n$, and by $D_l^*$ (resp. $D_l^u$) the number of balls of color $\mathfrak{d}$ in {\drawall} (resp. in the urn $u$) at time $l$.

We shall now prove that all the balls drawn out of {\drawall} after time $n$ are of color $\mathfrak{d}$ with a probability uniformly bounded away from 0.

Denote by $\mathcal{E}$ the set of urns that did not chose the color $\mathfrak{d}$, in symbols:
$$\mathcal{E}=\left\{u\in\mathcal{U}: \mathfrak{d}^u\neq\mathfrak{d}\right\}.$$
Then  $E=\left|\mathcal{E}\right|$ satisfies $(1-p)E<U/2$ (this is the hypothesis of the second case of the proof).

For all $l\geq n$ and $u\in\mathcal{U}$, set $k=lU+u$ and
$$Y_k=\mathbb{1}_{\{x_l^u=1\}\cup\{u\not\in\mathcal{E}\}}-\mathbb{1}_{\{x_l^u=-1\}\cap\{u\in\mathcal{E}\}}.$$

Set $k_1=(n+1)U$, and note that then the sequence $Y_{k_1+1}, Y_{k_1+2}, Y_{k_1+3},...$ satisfies the condition of Lemma \ref{lemme4} and is independent of $\mathcal{G}_1$. Set $\alpha=1/2-E(1-p)/U$ (any $\alpha\in(0,1-2E(1-p)/U)$ would suffice). Then we have
$$\mathbb{P}\left[A\right]=\theta_p^{\prime\prime}>0,$$
where $A=\{Y_{k_1+1}+ Y_{k_1+2}+...+ Y_{k_1+l}\geq\alpha l,~~\forall l\geq 0\}$.

For any $N\geq 1$ we have:

\begin{eqnarray*}
\mathbb{P}\left[\left.\bigcap_{l=n}^{n+N}\bigcap_{u\in\mathcal{U}}\left\{x_l^u=-1\right\}\cup\left\{\mathfrak{c}_l^u=\mathfrak{d}\right\}\right|\mathcal{F}_{n},\mathcal{G}_{-1}\right]&\geq&\mathbb{P}\left[\left.A\cap\bigcap_{l=n}^{n+N}\bigcap_{u\in\mathcal{U}}\left\{x_l^u=-1\right\}\cup\left\{\mathfrak{c}_l^u=\mathfrak{d}\right\}\right|\mathcal{F}_{n},\mathcal{G}_{-1}\right]\\
\end{eqnarray*}
\begin{eqnarray*}
&=&\mathbb{E}\left[\left.\mathbb{1}_A\times\prod_{l=n}^{n+N}\prod_{u\in\mathcal{U}}\mathbb{1}_{\left\{x_l^u=-1\right\}\cup\left\{\mathfrak{c}_l^u=\mathfrak{d}\right\}}\right|\mathcal{F}_{n},\mathcal{G}_{-1}\right]\\
&=&\mathbb{E}\left[\left.\mathbb{1}_A\times\prod_{l=n}^{n+N-1}\prod_{u\in\mathcal{U}}\mathbb{1}_{\left\{x_l^u=-1\right\}\cup\left\{\mathfrak{c}_l^u=\mathfrak{d}\right\}}\mathbb{E}\left[\left.\prod_{u\in\mathcal{U}}\mathbb{1}_{\left\{x_{n+N}^u=-1\right\}\cup\left\{\mathfrak{c}_{n+N}^u=\mathfrak{d}\right\}}\right|\mathcal{F}_{n+N-1},\mathcal{G}_{-1}\right]\right|\mathcal{F}_{n},\mathcal{G}_{-1}\right].\\
\end{eqnarray*}

Then on
$$A\cap\bigcap_{l=n}^{n+N-1}\bigcap_{u\in\mathcal{U}}\left\{x_l^u=-1\right\}\cup\left\{\mathfrak{c}_l^u=\mathfrak{d}\right\},$$
we have:

\begin{eqnarray*}
\mathbb{E}\left[\left.\prod_{u\in\mathcal{U}}\mathbb{1}_{\left\{x_{n+N}^u=-1\right\}\cup\left\{\mathfrak{c}_{n+N}^u=\mathfrak{d}\right\}}\right|\mathcal{F}_{n+N-1},\mathcal{G}_{-1}\right]&=&\prod_{u\in\mathcal{U}}\mathbb{E}\left[\left.\mathbb{1}_{\left\{x_{n+N}^u=-1\right\}\cup\left\{\mathfrak{c}_{n+N}^u=\mathfrak{d}\right\}}\right|\mathcal{F}_{n+N-1},\mathcal{G}_{-1}\right]\\
&=&\prod_{u\in\mathcal{U}}\mathbb{E}\left[\left.\mathbb{1}_{\left\{x_{n+N}^u=-1\right\}}+\mathbb{1}_{\left\{x_{n+N}^u=1\right\}\cap\left\{\mathfrak{c}_{n+N}^u=\mathfrak{d}\right\}}\right|\mathcal{F}_{n+N-1},\mathcal{G}_{-1}\right]\\
&=&\prod_{u\in\mathcal{U}}\mathbb{1}_{\left\{x_{n+N}^u=-1\right\}}+\mathbb{1}_{\left\{x_{n+N}^u=1\right\}}\mathbb{E}\left[\left.\mathbb{1}_{\left\{\mathfrak{c}_{n+N}^u=\mathfrak{d}\right\}}\right|\mathcal{F}_{n+N-1},\mathcal{G}_{-1}\right].\\
\end{eqnarray*}

Then on
$$\left(A\cap\bigcap_{l=n}^{n+N-1}\bigcap_{u\in\mathcal{U}}\left\{x_l^u=-1\right\}\cup\left\{\mathfrak{c}_l^u=\mathfrak{d}\right\}\right)\cap\left\{\mathfrak{c}_{n+N}^u=\mathfrak{d}\right\},$$
we have:

\begin{eqnarray*}
\mathbb{E}\left[\left.\mathbb{1}_{\left\{\mathfrak{c}_{n+N}^u=\mathfrak{d}\right\}}\right|\mathcal{F}_{n+N-1},\mathcal{G}_{-1}\right]&=&\frac{1}{1+\rho^{U(n+N)-2D_{n+N}^*}}\\
&\geq&\frac{1}{1+\rho^{-\alpha U N}}.\\
\end{eqnarray*}

Then:

\begin{eqnarray*}
\mathbb{E}\left[\left.\prod_{u\in\mathcal{U}}\mathbb{1}_{\left\{x_{n+N}^u=-1\right\}\cup\left\{\mathfrak{c}_{n+N}^u=\mathfrak{d}\right\}}\right|\mathcal{F}_{n+N-1},\mathcal{G}_{-1}\right]&\geq&\prod_{u\in\mathcal{U}}\left(\mathbb{1}_{\left\{x_{n+N}^u=-1\right\}}+\mathbb{1}_{\left\{x_{n+N}^u=1\right\}}\frac{1}{1+\rho^{-\alpha U N}}\right)\\
&\geq&\prod_{u\in\mathcal{U}}\left(\mathbb{1}_{\left\{x_{n+N}^u=-1\right\}}+\mathbb{1}_{\left\{x_{n+N}^u=1\right\}}\right)\frac{1}{1+\rho^{-\alpha U N}}\\
&\geq&\prod_{u\in\mathcal{U}}\frac{1}{1+\rho^{-\alpha U N}}\\
&=&\left(\frac{1}{1+\rho^{-\alpha U N}}\right)^U.\\
\end{eqnarray*}

Therefore by induction:

\begin{eqnarray*}
\mathbb{P}\left[\left.\bigcap_{l=n}^{n+N}\bigcap_{u\in\mathcal{U}}\left\{x_l^u=-1\right\}\cup\left\{\mathfrak{c}_l^u=\mathfrak{d}\right\}\right|\mathcal{F}_{n},\mathcal{G}_{-1}\right]&\geq&\theta_p^{\prime\prime}\prod_{l=0}^{N}\left(\frac{1}{1+\rho^{-\alpha U l}}\right)^U\\
&\geq&\theta_p^{\prime\prime}\prod_{l=0}^{\infty}\left(\frac{1}{1+\rho^{-\alpha U l}}\right)^U>0.\\
\end{eqnarray*}

Let $N$ go to infinity:

\begin{eqnarray*}
\mathbb{P}\left[\left.\bigcap_{l=n}^{\infty}\bigcap_{u\in\mathcal{U}}\left\{x_l^u=-1\right\}\cup\left\{\mathfrak{c}_l^u=\mathfrak{d}\right\}\right|\mathcal{F}_{n},\mathcal{G}_{-1}\right]&\geq&\theta_p^{\prime\prime}\prod_{l=0}^{\infty}\left(\frac{1}{1+\rho^{-\alpha U l}}\right)^U:=\eta_3>0.\\
\end{eqnarray*}

We can now conclude by the same argument that we used several times before. Set:
$$F_n=\bigcap_{l=n}^{\infty}\bigcap_{u\in\mathcal{U}}\left\{x_l^u=-1\right\}\cup\left\{\mathfrak{c}_l^u=\mathfrak{d}\right\},$$
and $F=\bigcup_{n=n_2}^\infty F_n$. We have:
$$\mathbb{1}_F=\lim_{n\in\infty}\mathbb{E}\left[\left.\mathbb{E}\left[\left.\mathbb{1}_{F}\right|\mathcal{F}_n, \mathcal{G}_{-1}\right]\right|\mathcal{F}_n\right]\geq\lim_{n\in\infty}\mathbb{E}\left[\left.\mathbb{E}\left[\left.\mathbb{1}_{F_n}\right|\mathcal{F}_n, \mathcal{G}_{-1}\right]\right|\mathcal{F}_n\right]\geq \eta_3>0.$$
Then $\mathbb{P}[F]=1$ which ends the proof of the second case.

Let us conclude this discussion with the case of equality $(1-p)B=U/2$. Since $W\leq B$ and $p>0$ we have:
$$U>(1-p)U = (1-p)B+(1-p)W = U/2 + (1-p)W$$
and then $(1-p)W<U/2$. Using the proof of the second case, it follows that for each time $n\geq n_2$ such that the majority color at time $n$ is black, there is a positive probability bounded away from $0$ that all the balls drawn out of all the $U$ urns combined after time $n$ are black. Therefore it suffice to prove that after $n_2$ there is infinitely many times $n$ such that black is the majority color, $\mathfrak{m}_n^*=\mathfrak{b}$. Now note that:
$$B_n^*-W_n^*=B_{n_2-1}^*-W_{n_2-1}^*+\sum_{l=n_2}^{n}\sum_{u\in\mathcal{U}}\mathbb{1}_{\{\mathfrak{c}^{u}_l=\mathfrak{b}\}}-\mathbb{1}_{\{\mathfrak{c}^{u}_l=\mathfrak{b}\}}\leq B_{n_2-1}^*-W_{n_2-1}^*+\sum_{l=n_2}^{n}\sum_{u\in\mathcal{U}}\left(\mathbb{1}_{A_l^u}-\mathbb{1}_{A_l^{uc}}\right),$$
where $A_l^u=\{X^{u}_l=-1\text{ and }\mathfrak{d}^u=\mathfrak{b}\}$. The last sum is a zero-drift random walk on $\Z$:
$$\mathbb{E}\left[\sum_{u\in\mathcal{U}}\mathbb{1}_{A_l^u}-\mathbb{1}_{A_l^{uc}}\right]=\mathbb{E}\left[-U+2\sum_{u\in\mathcal{U}}\mathbb{1}_{A_l^u}\right]=-U+2\times B \times \mathbb{P}\left[X_l^u=-1\right]=-U+2B(1-p)=0.$$
Therefore it becomes positive infinitely often. This proves that after a finite time, all the balls drawn out of all the $U$ urns combined are black.
\end{proof}

\subsection{Graphical summary}

We can now summarize the three cases in a single $(U,p)$ diagram. According to the Conformism Equation (\ref{CE}), when $p<1/2$, the interface between different phases are given by the straight lines $U=2(1-p)N_{\not c}$ with $N_{\not c}=1,~2,~3...$

\begin{center}
	\pspicture(0,-1)(13,9)

\psline[linewidth=1pt]{->}(0,0)(0,9)
\psline[linewidth=1pt]{->}(0,0)(13,0)

\psline[linewidth=.5pt](0,1)(6,.5)
\psline[linewidth=.5pt](0,2)(6,1)
\psline[linewidth=.5pt](0,3)(6,1.5)
\psline[linewidth=.5pt](0,4)(6,2)
\psline[linewidth=.5pt](0,5)(6,2.5)
\psline[linewidth=.5pt](0,6)(6,3)
\psline[linewidth=.5pt](0,7)(6,3.5)
\psline[linewidth=.5pt](0,8)(6,4)

\psline[linewidth=.9pt, linestyle=dotted](6,.5)(12,0)
\psline[linewidth=.9pt, linestyle=dotted](6,1)(12,0)
\psline[linewidth=.9pt, linestyle=dotted](6,1.5)(12,0)
\psline[linewidth=.9pt, linestyle=dotted](6,2)(12,0)
\psline[linewidth=.9pt, linestyle=dotted](6,2.5)(12,0)
\psline[linewidth=.9pt, linestyle=dotted](6,3)(12,0)
\psline[linewidth=.9pt, linestyle=dotted](6,3.5)(12,0)
\psline[linewidth=.9pt, linestyle=dotted](6,4)(12,0)

\psline[linewidth=1pt](6,0)(6,8)
\psline[linewidth=1pt, linestyle=dashed](0,3.5)(13,3.5)

\psline[linewidth=.5pt](-.09,.5)(.05,.5)
\psline[linewidth=.5pt](-.09,1)(.05,1)
\psline[linewidth=.5pt](-.09,1.5)(.05,1.5)
\psline[linewidth=.5pt](-.09,2)(.05,2)
\psline[linewidth=.5pt](-.09,2.5)(.05,2.5)
\psline[linewidth=.5pt](-.09,3)(.05,3)
\psline[linewidth=.5pt](-.09,3.5)(.05,3.5)
\psline[linewidth=.5pt](-.09,4)(.05,4)
\psline[linewidth=.5pt](-.09,4.5)(.05,4.5)
\psline[linewidth=.5pt](-.09,5)(.05,5)
\psline[linewidth=.5pt](-.09,5.5)(.05,5.5)
\psline[linewidth=.5pt](-.09,6)(.05,6)
\psline[linewidth=.5pt](-.09,6.5)(.05,6.5)
\psline[linewidth=.5pt](-.09,7)(.05,7)
\psline[linewidth=.5pt](-.09,7.5)(.05,7.5)
\psline[linewidth=.5pt](-.09,8)(.05,8)

\psline[linewidth=.5pt, linestyle=dotted](1.5,0)(1.5,3.5)
\psline[linewidth=.5pt, linestyle=dotted](3.6,0)(3.6,3.5)
\psline[linewidth=.5pt, linestyle=dotted](5,0)(5,3.5)

\rput(-.4,9){$U$}
\rput(13,-.4){$p$}
\rput(-.4,-.4){$0$}
\rput(6,-.4){$\frac{1}{2}$}
\rput(12,-.4){$1$}

\rput(1.5,-.4){$\frac{1}{8}$}
\rput(3.6,-.4){$\frac{3}{10}$}
\rput(5,-.4){$\frac{5}{12}$}

\rput(-.4,.5){$1$}
\rput(-.4,1){$2$}
\rput(-.4,1.5){$3$}
\rput(-.4,2){$4$}
\rput(-.4,2.5){$5$}
\rput(-.4,3){$6$}
\rput(-.4,3.5){$7$}
\rput(-.4,4){$8$}
\rput(-.4,4.5){$9$}
\rput(-.4,5){$10$}
\rput(-.4,5.5){$11$}
\rput(-.4,6){$12$}
\rput(-.4,6.5){$13$}
\rput(-.4,7){$14$}
\rput(-.4,7.5){$15$}
\rput(-.4,8){$16$}

\rput(2.7,.4){$N_{\not c}=0$}
\rput(2.7,1.1){$N_{\not c}\leq 1$}
\rput(2.7,1.9){$N_{\not c}\leq 2$}
\rput(2.7,2.7){$N_{\not c}\leq 3$}
\rput(2.1,3.75){$N_{\not c}\leq 4$}
\rput(2.1,4.5){$N_{\not c}\leq 5$}
\rput(2.1,5.3){$N_{\not c}\leq 6$}
\rput(2.1,6.1){$N_{\not c}\leq 7$}

\rput(9,2.9){$N_{\not c}=0$}
	\endpspicture
	
\noindent
{\footnotesize Figure 7:
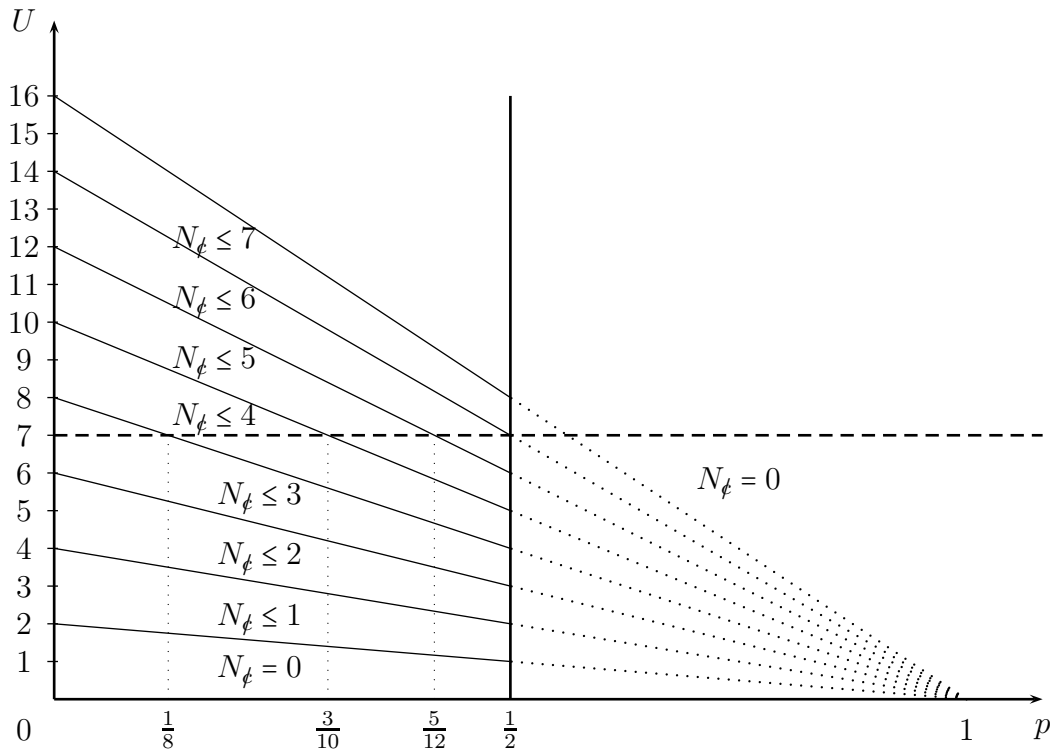
Summary of different cases depending on $p$ and $U$.}
\setcounter{figure}{7}

\vspace{0.3cm}
\noindent
	\end{center}

Take for example $U=7$ (cf. the dashed horizontal line) then
\begin{itemize}
\item if $p\in[0,1/8]$ there is at most three nonconformist urns (and there is actually a positive probability that $N_{\not c}=i$ for each $i=0,1,2$ or $3$);
\item if $p\in(1/8,3/10]$ there is at most four nonconformist urns;
\item if $p\in(3/10,5/12]$ there is at most five nonconformist urns;
\item if $p\in(5/12,1/2)$ there is at most six nonconformist urns;
\item if $p\in[1/2,1]$ there is no nonconformist urns.
\end{itemize}

\section{Two generalizations}

\subsection{Extension to more than two colors}

Let us now consider what happens if there are more than two colors available. Define $\mathfrak{C}=(\mathfrak{d}_1,\mathfrak{d}_2,\dots,\mathfrak{d}_C)$ the set of $C$ possible colors, with $C\in\N$. The Interacting Urns Mechanism with $C\geq 2$ colors is just the natural generalization of the $C=2$ case. The probability to draw a ball of color $\mathfrak{d}_i$ from an urn that contains $N_j$ balls of color $\mathfrak{d}_j$ for $1\leq j\leq C$, equals:
$$\frac{w_{N_i}}{\sum_{i=1}^Cw_{N_j}}.$$
\textbf{Claim:} Theorem \ref{psur} (supercritical case), Theorem \ref{pcritic} (critical case), Proposition \ref{propsub} and Theorem \ref{thsub} (subcritical case) are still true when $C\geq 2$.\\
If $p<1/2$ an $X$-suitable {\up} in the case $C\geq 2$ can be schematically represented as follows:
	\begin{center}
	\pspicture(0,-1)(8.5,2.5)

			\psline[border=0pt](0,2)(0,0)
			\psline[border=0pt](0,0)(1,0)
			\psline[border=0pt](1,0)(1,2)
			\psline[border=0pt](0,.6)(1,.6)
			\psline[border=0pt](0,1.8)(1,1.8)
			
			\psline[border=0pt]{-}(1.5,2)(1.5,0)
			\psline[border=0pt]{-}(1.5,0)(2.5,0)
			\psline[border=0pt]{-}(2.5,0)(2.5,2)
			\psline[border=0pt](1.5,.6)(2.5,.6)
			\psline[border=0pt](1.5,1.8)(2.5,1.8)

			\psframe[linecolor=black,fillcolor=gray,fillstyle=solid](3,.6)(4,1.8)
			\psline[border=0pt]{-}(3,2)(3,0)
			\psline[border=0pt]{-}(3,0)(4,0)
			\psline[border=0pt]{-}(4,0)(4,2)
			
			\psframe[linecolor=black,fillcolor=black,fillstyle=solid](4.5,.6)(5.5,1.8)
			\psline[border=0pt]{-}(4.5,2)(4.5,0)
			\psline[border=0pt]{-}(4.5,0)(5.5,0)
			\psline[border=0pt]{-}(5.5,0)(5.5,2)

			\psframe[linecolor=black,fillcolor=gray,fillstyle=solid](6,.6)(7,1.8)
			\psline[border=0pt]{-}(6,2)(6,0)
			\psline[border=0pt]{-}(6,0)(7,0)
			\psline[border=0pt]{-}(7,0)(7,2)
			
			\psframe[linecolor=black,fillcolor=lightgray,fillstyle=solid](7.5,.6)(8.5,1.8)
			\psline[border=0pt]{-}(7.5,2)(7.5,0)
			\psline[border=0pt]{-}(7.5,0)(8.5,0)
			\psline[border=0pt]{-}(8.5,0)(8.5,2)
			
			\rput(1.25,-.2){$\underbrace{~~~~~~~~~~~~~~~~~~}$}
			\rput(5.75,-.2){$\underbrace{~~~~~~~~~~~~~~~~~~~~~~~~~~~~~~~~~~~~~~~~~}$}
			\rput(1.25,-.5){{conformist urns}}
			\rput(5.75,-.5){{nonconformist urns}}
			\rput(-.1,1.3){$\left\{\raisebox{0pt}[20pt][0pt]{~}\right.$}
			\rput(-.1,.3){$\left\{\raisebox{0pt}[7.5pt][0pt]{~}\right.$}
			\rput(-.8,1.2){$1-p$}
			\rput(-.5,.3){$p$}

	\endpspicture
	
\noindent
{\footnotesize Figure 8:
Schematic representation of an $X$-suitable {\up} when $C\geq 2$.}
\setcounter{figure}{8}

\vspace{0.3cm}
\noindent
	\end{center}

Let us justify briefly the claim and then explain how to adapt the proofs so that they work with $C\geq 2$.

In the case of two colors we used arguments of the following form: ``Consider the majority color at time $n$, and prove that there is a positive probability uniformly bounded away from $0$ that after time $n$ we draw only this color''. But let us think about what the term {\em``majority color''} means in our particular settings. In politics, for instance in a parliament, a majority could be {\em absolute} or {\em relative}. In the first case, it means ``having more than $50\%$ of the seats'' whereas in the second it only means ``having more seats than any other group''. In the context of our proof, the relevant interpretation is the second: what matters is that a certain color has more balls than any other color.

For example, consider the two following configurations for a single urn with two colors:
\begin{enumerate}[i)]
\item $n=10$ and $(B_{10},W_{10})=(9,1)$ ;
\item $n=1000$ and $(B_{1000},W_{1000})=(505,495)$.
\end{enumerate}
Then the probability for the black to win is larger in the second case ($B-W=10$) than in the first one ($B-W=8$). The information that the proportions are $(90\%,10\%)$ in the first case and $(50.5\%,49.5\%)$ in the second is not relevant.

This finding is important since, for $C\geq 3$, there is not necessarily a color with absolute majority, but there will always be at least one color with relative majority. Let us take for example the case of one single urn with three colors, and suppose without loss of generality that the majority color is $\mathfrak{d}_1$, that is, $N_1\geq N_2$ and $N_1\geq N_3$. Then the probability to draw only balls of color $\mathfrak{d}_1$ starting from this configuration is:

$$P=\prod_{i=0}^{\infty}\frac{\rho^{N_1+i}}{\rho^{N_1+i}+\rho^{N_2}+\rho^{N_3}}=\prod_{i=0}^{\infty}\frac{1}{1+\rho^{N_2-N_1-i}+\rho^{N_3-N_1-i}}\geq\prod_{i=0}^{\infty}\frac{1}{1+\rho^{-i}+\rho^{-i}}=\prod_{i=0}^{\infty}\frac{1}{1+2\rho^{-i}}>0.$$
If there were $C$ colors instead, we would have had:
$$P\geq\prod_{i=0}^{\infty}\frac{1}{1+(C-1)\rho^{-i}}>0$$

In other words, all our draw-the-majority-color-with-positive-probability arguments in the last section still work and all the results are still true.

\subsection{Extensions to a wider class of reinforcement weight sequences}

Let us consider the following hypothesis about the reinforcement weight sequence:
\begin{equation}
\label{H1}
\tag{$H_1$}
\liminf_{i\to\infty}\frac{w_{i+1}}{w_i}>1.
\end{equation}

\textbf{Claim:} Theorem \ref{psur} (supercritical case), Theorem \ref{pcritic} (critical case), Proposition \ref{propsub} and Theorem \ref{thsub} (subcritical case) are still true under (\ref{H1}).\\

We can immediately make two remarks about this hypothesis:
\begin{enumerate}[i)]
\item The exponential reinforcement with $\rho>1$ satisfies (\ref{H1}):
$$\liminf_{i\to\infty}\frac{\rho^{i+1}}{\rho^i}=\liminf_{i\to\infty}\rho=\rho>1$$
\item Any reinforcement weight sequence that satisfies (\ref{H1}) satisfies (\ref{SRP}).
\end{enumerate}

Under (\ref{H1}), and for $a,b\geq I$ and $a>b$, we get a lower bound for the quantity $\pi(a,b)$:
$$\pi(a,b)=\frac{w_a}{w_a+w_b}=\frac{1}{1+\frac{w_b}{w_a}}\geq\frac{1}{1+\frac{w_b}{w_b\rho^{a-b}}}=\frac{1}{1+\rho^{b-a}}.$$
The last quantity corresponds to the exponential reinforcement. This inequality can be used in the arguments each time we need a lower bound on $\pi(a,b)$, and allows us to finish the proofs just as in the exponential case. Note that this lower bound is true only if $a>b$, but it is sufficient since we only use $\pi(a,b)$ when $a$ is already the number of balls of the majority color.

\section{Conclusion}

A natural follow-up to the results of this paper is to search for the distribution of the number of nonconformist urns. That is, answering the question ``What is the value of
$$\mathbb{P}\left[N_{\not c}=i\right]$$
for $i<U/2(1-p)$?'' This seems to be a difficult question. Another natural question is to wonder what happens if the reinforcement sequence is strongly  but not exponentially reinforced, for example if it grows only polynomially.

We can actually conjecture that under (\ref{SRP}) a Reinforced Urn Process can behave in two different ways: either it behaves as the exponentially reinforced urn processes (which means that Theorem \ref{psur} (supercritical case), Theorem \ref{pcritic} (critical case), Proposition \ref{propsub} and Theorem \ref{thsub} (subcritical case) are true) or all the urns draw the same color after a finite time, whenever $p>0$. Roughly speaking, either the reinforcement sequence is strong enough to allow individual urns to have some independence for being nonconformist, or all the urns eventually have to give in to the majority color. It also seems natural to conjecture that polynomial reinforcements behave in the second way. It is actually easy to see that they could not behave as the exponential reinforcements.

It turns out that already proving that all the urns eventually draw the same color in the case $p=1$ and $w$ general, that is when the $U$ urns behave just like one big urn with $U$ draws at each step, is not entirely trivial. This work is in progress and shall be the subject of a forthcoming paper.
\\

When we study reinforcement processes it is natural to wonder if there is a natural way to say that one reinforcement sequence is stronger than an other one. In other words, is there a natural way to put an order $\leq$ over the set of reinforcement sequences? The (\ref{SRP}) gives a first rough answer dividing this set in two, depending on whether the sum of inverses is finite or not. The results of this paper allows us to split again the set of strongly reinforced sequences in two parts: those that satisfy Theorem \ref{psur}, Theorem \ref{pcritic}, Proposition \ref{propsub} and Theorem \ref{thsub} and those that do not. It would be interesting to know if this splitting could be expressed in a natural way independent of Interacting Urns Models as (\ref{SRP}) do. And could this splitting have a meaning in other reinforced processes as for example, in the context of Reinforced Random Walks?

\textbf{Acknowledgments.} I would like to thank my Ph.D. advisor, Vlada Limic, for many insightful discussions about this work as well as detailed comments on earlier drafts.

\bibliographystyle{plain}
\bibliography{biblio}

\begin{thebibliography}{1}

\bibitem{davis}
B.~Davis.
\newblock Reinforced random walk.
\newblock {\em Prob. Th Rel. Fields}, 84:203--229, 1990.

\bibitem{durrett}
R.~Durrett.
\newblock {\em Probability: Theory and Examples}.
\newblock Cambridge U. Press, 4th edition, 2010.

\bibitem{sabot1}
N.~Enriquez and C.~Sabot.
\newblock Edge oriented reinforced random walks and rwre.
\newblock {\em C. R. Math. Acad. Sci.}, 335(11):941--946, 2002.

\bibitem{sabot2}
N.~Enriquez and C.~Sabot.
\newblock Random walks in a dirichlet environment.
\newblock {\em Electron. J. Probab.}, 11(31):802--817, 2006.

\bibitem{limictarres}
V.~Limic and P.~Tarr\`{e}s.
\newblock Attracting edge and strongly edge reinforced random walks.
\newblock {\em Ann. Probab.}, 35(5):1783--1806, 2007.

\bibitem{merkl}
F.~Merkl and S.~W.~W. Rolles.
\newblock Recurrence of edge-reinforced random walk on a two-dimensional graph.
\newblock {\em Ann. Probab.}, 37(5):1679--1714, 2009.

\bibitem{pemantle_survey}
R.~Pemantle.
\newblock A survey of random processes with reinforcement.
\newblock {\em Probability surveys}, 4:1--79, 2007.

\bibitem{polya}
G.~P\'{o}lya.
\newblock Sur quelques points de la th\'{e}orie des probabilit\'{e}s.
\newblock {\em Ann. Inst. H. Poincar\'{e}}, 1:117--161, 1930.

\bibitem{toth}
B.~T\'{o}th.
\newblock Limit theorems for weakly reinforced random walks.
\newblock {\em Studia Sci. Math. Hungar.}, 33:321--337, 1997.

\end{thebibliography}

\end{document}